\documentclass [a4paper, 12pt] {amsart}

\usepackage {hyperref}

\usepackage [utf8] {inputenc}

\usepackage{mathtext}

\usepackage {version}

\usepackage {amsmath}
\usepackage {amssymb}
\usepackage {amsthm}

\usepackage {graphicx}

\usepackage {appendix}

\usepackage {url}

\usepackage{mathtext}
\usepackage{amsmath}
\usepackage{amsfonts}
\usepackage{amssymb}
\usepackage{mathrsfs}
\usepackage{amsthm}

\usepackage {needspace}

\excludeversion {svntags}

\newtheorem {theorem} {Theorem}

\newtheorem* {theoremnn} {Theorem}

\newtheorem {lemma} [theorem] {Lemma}
\newtheorem {proposition} [theorem] {Proposition}% [section]
\newtheorem {definition} [theorem] {Definition}% [section]
\newtheorem {corollary} [theorem] {Corollary}% [section]

\newcommand {\apclass} [1] {\ensuremath{\mathrm A_{#1}}}

\newcommand {\lclass} [2] {\ensuremath{\mathrm L_{#1} \left( #2 \right) }}
\newcommand {\lsclass} [1] {\ensuremath{\mathit l^{#1} }}

\newcommand {\lclassg} [1] {\ensuremath{\mathrm L_{#1}}}

\newcommand {\hclassg} [1] {\ensuremath{\mathrm H_{#1}}}

\newcommand {\BMO} {\ensuremath {\mathrm {BMO}}}

\DeclareMathOperator* {\essinf} {ess\,inf}
\DeclareMathOperator* {\esssup} {ess\,sup}
\DeclareMathOperator* {\supp} {supp}

\excludeversion {comment}

\newcommand {\weightu} {\ensuremath {\mathit u}}
\newcommand {\weightv} {\ensuremath {\mathit v}}
\newcommand {\weightw} {\ensuremath {\mathit w}}

\newcommand {\xalphab} {\ensuremath {{X \strut}^\alpha {\lclassg {1}\hspace{-0.18cm}\strut}^{1 - \alpha}}}

\newcommand {\abp} [2] {\ensuremath {B\apclass {#1} \left(#2\right)}}

\excludeversion {abridged}
\includeversion {unabridged}

\begin {document}

\title [$\apclass {1}$-regularity] {$\mathbf A_1$-regularity and boundedness\\of Calderon-Zygmund operators}
\author {D.~V.~Rutsky}
\email {rutsky@pdmi.ras.ru}
\date {\today}
\address {St.Petersburg Department
of Steklov Mathematical Institute RAS
27, Fontanka
191023 St.Petersburg
Russia}

\keywords {$\apclass {p}$-regularity, Calderon-Zygmund operator, Hardy-Littlewood maximal operator}
%\keywords {$\apclass {p}$-regularity, Calderon-Zygmund operators, dyadic averages operator, harmonic conditional expectation operator}

\begin {abstract}
The Coifman-Fefferman inequality implies quite easily
that a Calde\-ron-Zygmund operator $T$ acts boundedly in a Banach lattice $X$
on $\mathbb R^n$ if the Hardy-Littlewood maximal operator $M$ is bounded in both $X$ and $X'$.
\begin {abridged}
We 
\end {abridged}
\begin {unabridged}
 %In this paper
We discuss this phenomenon in some detail
and 
\end {unabridged}
establish a converse result under the assumption that $X$ satisfies the Fatou property and $X$
is $p$-convex and $q$-concave with some $1 < p, q < \infty$:
if a linear operator $T$ is bounded in $X$ and $T$ is nondegenerate in a certain sense (for example, if $T$ is a Riesz transform)
then $M$ is bounded in both $X$ and $X'$.
\end {abstract}

\maketitle

\setcounter {section} {-1}

\section {Introduction}

\label {anintroduction}

\begin {unabridged}
The problem of characterizing the spaces in which (and between which) the operators of harmonic analysis act boundedly %is, perhaps,
 %the central one
lies in the core of the modern harmonic analysis, and it definitely has far-reaching consequences in terms of applications.
These operators in a vast number of cases can be represented by (or the corresponding questions reduced to the study of)
a general Calderon-Zygmund operator.  The study of such operators has received a lot of attention over the past several decades and
significant advancements have been made.
To mention a few highlights: the quest for practical conditions that guarantee boundedness of a Calderon-Zygmund operator in $\lclassg {2}$ led to
useful $\mathrm T 1$ theorems, new approaches to the classical proofs %allowed
have made it possible
to significantly relax the doubling condition on the underlying
measurable space, the action of such operators was studied in detail in a wide variety of spaces beyond the classical Lebesgue spaces $\lclassg {p}$,
and a number of representations for such operators were developed together with highly refined techniques that recently
 %allowed
%to solve
yielded answers to
several long-standing %open
problems such as the $\apclass {2}$-hypothesis (positive) and the $\apclass {1}$ conjecture of
Muckenhoupt and Wheeden (negative).  Although it seems that the focus has always been on particular classes of spaces, weighted Lebesgue spaces
$\lclass {p} {\weightw}$ being of a particular interest (not least because of their rather general nature which has long been noted),
results extending various useful relationships to fairly general classes of spaces, and indeed sometimes demonstrating exhaustively the true
scope of what has been known for many years, recently began to emerge.
\end {unabridged}

The purpose of the present work is to establish the following theorem showing that
the boundedness of Calderon-Zygmund singular integral operators $T$ and
the boundedness of the Hardy-Littlewood maximal operator $M$ in both the lattice and its dual
is actually the same property in a fairly general class of Banach lattices.
It constitutes a substantial improvement over the respective results of \cite {rutsky2011en}.

The (standard) definitions and basic facts concerning Banach lattices and Calderon-Zygmund operators can be found in Section~\ref {preliminaries}.
The notion of an $\apclass {2}$-nondegenerate lattice is introduced in Defition~\ref {nondego} below;
for now we say that $R$ can be any of the Riesz transforms $\{R_j\}_{j = 1}^n$.
We fix a $\sigma$-finite measurable space $\left(\Omega, \mu\right)$
which we understand as a space for the second variable $\omega$ in $(x, \omega) \in \mathbb R^n \times \Omega$
(unless indicated otherwise, all operators are assumed to act in the first variable $x$ only);
this allows us to naturally include lattices with mixed norm such as $X (\lsclass {r})$ in this setting.
 %where $X$ is a lattice on $\mathbb R^n$.
\begin {theorem}
\label {themcr}
Suppose that $X$ is a Banach lattice of measurable functions on $\mathbb R^n \times \Omega$ that satisfies the Fatou property
and $X$ is $p$-convex and $q$-concave with some $1 < p, q < \infty$.  Let $R$ be a Calderon-Zygmund operator in $\lclass {2} {\mathbb R^n}$ such that
both $R$ and $R^*$ are $\apclass {2}$-nondegenerate.
The following conditions are equivalent.
\begin {enumerate}
\item
The Hardy-Littlewood maximal operator $M$ acts boundedly in $X$ and in the order dual $X'$ of $X$.
\item
All Calderon-Zygmund operators act boundedly in $X$.
\item
$R$ acts boundedly in $X$.
\end {enumerate}
\end {theorem}

\begin {abridged}
Implication $1 \Rightarrow 2$ is established by Proposition~\ref {czobounda1p} in Section~\ref {preliminaries}.
Although it is hard to come by this sufficient condition for boundedness of Calderon-Zygmund operators in the literature,
it is certainly not new; see \cite [Remark~4.3] {karlovichlerner2005}.
Implication $2 \Rightarrow 3$ is trivial.
Implication $3 \Rightarrow 1$ is established by Theorem~\ref {a1necc} in Section~\ref {noaro}.
 %Although
The argument itself is technically rather simple; however,
it relies heavily on the theory of $\apclass {p}$-regular Banach lattices, a part of which we develop further in Section~\ref {aprmainlemma},
and the proof taken as a whole involves overall two distinct applications of the Ky-Fan--Kakutani fixed point theorem and a variant of
the Maurey--Krivine factorization theorem which is based on the Grothendieck theorem.
\end {abridged}

\begin {unabridged}
We will explore several proofs of implication~$1 \Rightarrow 2$ in Section~\ref {czos} below.
Although it is hard to come by this sufficient condition for boundedness of Calderon-Zygmund operators in the literature,
it is certainly not new; see, e.~g., \cite [Remark~4.3] {karlovichlerner2005}.
Implication~$2 \Rightarrow 3$ is trivial, and implication~$3 \Rightarrow 1$, which is in a sense the main point of the present work,
is established in Section~\ref {noaro};
although the argument itself is technically simple,
it relies heavily on the theory of $\apclass {p}$-regular Banach lattices, a part of which we develop further in Section~\ref {aprmainlemma},
and the proof taken as a whole involves overall two distinct applications of the Ky-Fan--Kakutani fixed point theorem and a variant of
the Maurey--Krivine factorization theorem which is based on the Grothendieck theorem\footnote {
A note of caution concerning how this paper is laid out seems to be necessary:
since the author tried to properly introduce and discuss at length all elements (more or less well-known with possibly a few exceptions)
leading to this result in order to explore
possible connections and extensions, it was convenient to postpone the main argument until the very end.
Thus the impatient reader who wishes to study the proof of the implication $3 \Rightarrow 1$ is advised to skip right away
down to Section~\ref {aprmainlemma}
and refer to the rest of the paper and to \cite {rutsky2011en} as necessary.
An abridged version of this paper made for submission to a journal is also available upon request (or by configuring the \LaTeX{} sources in a certain way).
}.

As it will be seen,
the sufficiency of Condition~3 of Theorem~\ref {themcr} for the other conditions
actually extends to a wide class of singular operators that are nondegenerate in a certain sense.
 %Although for simplicity we present the main result for $\mathbb R^n$,
%it is easily seen that Theorem~\ref {themcr} is valid for spaces of homogeneous type.
The proof, which is covered by Proposition~\ref {czobounda1p} in Section~\ref {czos} and by
Theorem~\ref {a1necc} in Section~\ref {noaro} below can easily be generalized to the case of a general space of homogeneous type
instead of just $\mathbb R^n$ if there exists a suitable nondegenerate operator $R$; %that can play the part of the Riesz transform $R_j$.
it is not clear whether every space of homogeneous type has at least one such operator.
It is easy to see that the proof of Theorem~\ref {themcr} also works in the vector-valued case,
i.~e. for lattices of measurable functions like $X (\lsclass {r})$, where $X$ is a lattice on $\mathbb R^n$.
The $p$-convexity and $q$-concavity assumptions are probably not necessary (they are not used in the implication $1 \Rightarrow 2$)
and
I conjecture that they in themselves are a consequence of any of the conditions
of Theorem~\ref {themcr}; that Condition~1 implies $p$-convexity and $q$-concavity with some $1 < p, q < \infty$
is known to hold true at least in the case of the variable exponent Lebesgue spaces (see, e.~g., \cite [Theorem~4.7.1] {varpbook}),
and it seems that it is possible to adapt the same argument to cover suitable nondegenerate singular integral operators as well.
Recently in \cite [Theorem~5.42] {varlsp} it was established that if all Riesz transforms $R_j$ are bounded in $\lclassg {p (\cdot)}$ then
the exponent $p (\cdot)$ is bounded away from $1$ and $\infty$ (and thus lattice $\lclassg {p (\cdot)}$ satisfies the $p$-convexity and $q$-concavity
assumptions in this case).
It is also interesting to note that implication $1 \Rightarrow 2$ easily extends in a certain natural way
to the case of operators acting between different Banach lattices; see Theorem~\ref {dpmxy} in Section~\ref {czos} below.
\end {unabridged}

\begin {abridged}
% FIXME
\end {abridged}

\begin {unabridged}
Let us briefly outline some the contributions that led to this result.  In the standard part of the theory describing
the properties of the Calderon-Zygmund operators
(see, e.~g., \cite {stein1993}) in the Lebesgue space setting the maximal operator plays an essential part.
In the case of weighted Lebesgue spaces $\lclass {p} {\weightw}$
Theorem~\ref {themcr}, of course, follows from the theory of the Muckenhoupt weights (see, e.~g., \cite [Chapter~5] {stein1993}) that
individually links the conditions
of Theorem~\ref {themcr} to the Muckenhoupt condition of the weight $\weightw$.
Of a particular interest in this regard is the Coifman-Fefferman inequality \cite {coifmanfefferman1974}
\begin {equation}
\label {coifmani}
\int |T f|^p \omega \leqslant C \int (M f)^p \omega, \quad 0 < p < \infty,
\end {equation}
with $C$ independent of $f$,
that holds true for Calderon-Zygmund operators $T$ %(with the standard conditions on the kernel)
and any weight $\omega \in \apclass {\infty}$
for all locally summable functions $f$ such that the right-hand part of \eqref {coifmani} is finite.
Thus $T$ is estimated in terms of $M$ for a relatively
wide class of Muckenhoupt weights even though $M$ may not act boundedly
in the corresponding weighted Lebesgue space.  There is a large number of various extensions of \eqref {coifmani}; see, e.~g.,
\cite {weightsandextrapolation}.  On the other hand, %application of
making use of the
duality and the famous construction due to Rubio de Francia
allowed a large number of very useful extrapolation results
that essentially exploit a very simple idea: if $M$ is bounded in $X$ then any $f \in X$ can be pointwise dominated with a controlled increase of norm
by some weight $\weightw \in \apclass {1}$, and the converse is also trivially true.
It is natural to call such lattices \emph {$\apclass {1}$-regular} by analogy with $\BMO$-regularity (see \cite {rutsky2011en}).
For example, this idea works very well in the case of variable exponent Lebesgue spaces $\lclassg {p (\cdot)}$ where
the behavior of boundedness of $M$ under duality and certain scaling operations is nice and well understood; see, e.~g., \cite [\S7.2] {varpbook}.
The Coifman-Fefferman inequality~\eqref {coifmani} with $p = 1$
gives a very easy proof of the implication $1 \Rightarrow 2$ of Theorem~\ref {themcr}; see Proposition~\ref {czobounda1p} in Section~\ref {czos} below.
And this is far from the only way to establish this implication;
we will also discuss in Section~\ref {czos} below how
some of the recent results by A.~Lerner \cite {lerner2004}, \cite {lerner2010} and \cite {lerner2012}
 also give the necessary tools
to effortlessly establish this implication.

The study of the duality of $\BMO$-regularity, initially motivated by certain problems in the theory of interpolation of Hardy-type spaces,
eventually led in \cite {rutsky2011en}
to a refinement and generalization to the general spaces of homogeneous type
of certain properties and results concerning the interplay of various majorization and boundedness properties
that were previously known only in the case of the unit circle~$\mathbb T$.
In particular, the main result of \cite {rutsky2011en} is similar to Theorem~\ref {themcr} because it links boundedness of $T$ and $M$
in lattices of the form
$\xalphab$
for $0 < \beta < 1$ and sufficiently small $0 < \alpha < 1$ to another property,
namely to $\BMO$-regularity of $X$.
This, admittedly, still left much more to be desired in terms of refinements,
 %because
since
unlike $\apclass {1}$-regularity the $\BMO$-regularity property,
which proved to be very useful in certain questions pertaining to spaces on the unit circle $\mathbb T$,
so far does not seem to be as useful in the case of the spaces on $\mathbb R^n$ in the same capacity.
The results of the present work can be regarded as an extension and an application
of the techniques described in \cite {rutsky2011en}.
 %In this paper we will see how the techniques described in \cite {rutsky2011en} can be adapted to establish the converse implications
%of Theorem~\ref {themcr}.
\end {unabridged}

The paper is organized as follows.
\begin {abridged}
In Section~\ref {preliminaries} we provide the definitions and basic properties that will be used in the text
and prove the implication $1 \Rightarrow 2$ of Theorem~\ref {themcr}.
Section~\ref {aprmainlemma} contains a new sufficient condition for $\apclass {1}$-regularity.
In Section~\ref {noaro} we prove the converse implication $3 \Rightarrow 1$ of Theorem~\ref {themcr}.
Further remarks are given in Section~\ref {concrems}.
\end {abridged}
\begin {unabridged}
In Section~\ref {preliminaries} we introduce some basic notions pertaining to Banach lattices and
spaces of homogeneous type.
In Section~\ref {mwaam} certain known facts about Muckenhoupt weights and $\apclass {p}$-regular spaces are outlined.
In Section~\ref {czos} we briefly describe Calderon-Zygmund operators and show
several ways
to obtain the implication $1 \Rightarrow 2$ of Theorem~\ref {themcr}.
Then in Section~\ref {nsos} we discuss some results having to do with operators that are nondegenerate in a certain sense.
Section~\ref {aprmainlemma} contains a new result that gives a sufficient condition for a lattice $X$ to be $\apclass {1}$-regular
in terms of $\apclass {1}$-regularity of lattice $X^\delta$ and $\apclass {p}$-regularity of lattice $X$.
Finally, in Section~\ref {noaro} we prove the converse implication $3 \Rightarrow 1$ of Theorem~\ref {themcr}.
\end {unabridged}

\section {Preliminaries}

\label {preliminaries}

%\begin {abridged}
In this section we briefly go over the basic definitions and facts used by the rest of the paper.
For the generalities on
\begin {abridged}
the real harmonic analysis see, e.~g., \cite {denghan}, \cite {stein1993};
for
\end {abridged}
the Banach lattices and their
properties see, e.~g., \cite [Chapter~10] {kantorovichold},
\cite {cbs}.
A space of homogeneous type $(S, \nu)$ is a 
quasimetric space equipped with a Borel measure $\nu$ that has the doubling property, i.~e.
$\nu (B (x, 2 r)) \leqslant c \nu (B (x, r))$
for all $x \in S$ and $0 < r < \infty$ with some constant $c$, where $B (x, r)$ is the ball of radius $r$ centered at $x$.
The main example here is $S = \mathbb R^n$
equipped with the Lebesgue measure.

A quasi-normed lattice of measurable functions $X$ is a quasi-normed space of measurable functions $X$
in which the norm is compatible with the natural order; that is, if $|f| \leqslant g$ a.~e. for some function $g \in X$
then $f \in X$ and $\|f\|_X \leqslant \|g\|_X$.
For simplicity we only work with lattices $X$ such that $\supp X = S \times \Omega$.
For a Banach lattice of measurable functions $X$, any order continuous functional $f$ on $X$
(order continuity is understood in the sense that for any sequence $x_n \in X$
such that $\sup_n |x_n| \in X$ and $x_n \to 0$ a. e. one also has $f (x_n) \to 0$) has an integral representation
$f (x) = \int x y_f$ for some measurable function $y_f$ which can be identified with $f$.
The set of all such functionals $X'$ is a Banach lattice with
the norm defined by $\|f\|_{X'} = \sup_{g \in X, \|g\|_X = 1} \int |f g|$.  The lattice $X'$ is called the order dual of the lattice $X$.
The norm of a lattice $X$ is said to be order continuous if for any nonincreasing sequence
$x_n \in X$ converging to $0$ a. e. one also has
$\|x_n\|_X \to 0$.  Order continuity of the norm of a Banach lattice $X$
is equivalent to $X^* = X'$, and it is also equivalent to density of the simple functions in $X$.
A lattice $X$ has the Fatou property if for any
$f_n, f \in X$ such that $\|f_n\|_X \leqslant 1$ and the sequence $f_n$ 
converges to $f$ a. e. it is also true that $f \in X$ and $\|f\|_X \leqslant 1$.
The Fatou property of a lattice $X$ is equivalent to $(\nu \times \mu)$-closedness of the unit ball
$B_X$ of the lattice $X$
(here and elsewhere $(\nu \times \mu)$-convergence denotes the
convergence in measure in any measurable set $E$ such that $(\nu \times \mu) (E) < \infty$).
If the lattice $X$ is Banach then the Fatou property is equivalent to order reflexivity of $X$, i.~e. to the relation $X'' = X$.
For a lattice $X$ either one of the Fatou property or the order continuity of norm property
is sufficient to guarantee that the lattice $X'$ is norming for
$X$, i.~e. that $\|f\|_X = \sup_{g \in X', \|g\|_{X'} = 1} \int f g$ for all $f \in X$.

For any two quasi-normed lattices $X$ and $Y$ on the same measurable space
the set of poinwise products of their functions
$$
X Y = \{ f g \mid f \in X,\break g \in Y\}
$$
is a quasi-normed lattice with the norm defined by
$$
\|h\|_{X Y} = \inf_{h = f g} \|f\|_X \|g\|_Y.
$$
If both lattices $X$ and $Y$ satisfy the Fatou property then the lattice $X Y$ also has the Fatou property.
If either of the lattices $X$ and $Y$ has order continuous quasi-norm then the quasi-norm of $X Y$ is also order continuous.

For any $\delta > 0$ and a quasi-normed lattice $X$ the lattice $X^\delta$ consists of all measurable functions $f$
such that $|f|^{1/\delta} \in X$ with a quasi-norm $\|f\|_{X^\delta} = \| |f|^{1/\delta} \|_X^\delta$.
For example, $\lclassg {p}^\delta = \lclassg {\frac p \delta}$.
It is easy to see that $(X Y)^\delta = X^\delta Y^\delta$ for any $X$, $Y$ and $\delta$,
and $X^\delta$ naturally inherits many properties from $X$.
For any $0 < \delta \leqslant 1$, if $X$ is a Banach lattice then
$X^\delta$ is also a Banach lattice.
If $X$ and $Y$ are Banach lattices then for any $0 < \delta < 1$ lattice
$X^{1 - \delta} Y^\delta$, sometimes called the \emph {Calderon-Lozanovsky product} of $X$ and $Y$, is also Banach; moreover, one has a very useful relation
$(X^{1 - \delta} Y^\delta)' = (X')^{1 - \delta} (Y')^\delta$ (see \cite {calderon1964}, \cite {lozanovsky1969}).
If $Z = X^{1 - \delta} Y^\delta$ has either the Fatou property or order continuous norm
then
$Z$ is an exact interpolation space of exponent $\delta$ between $X$ and $Y$;
see, e.~g., \cite {lozanovsky1972},~\cite {calderon1964},~\cite {kps}.

Let $1 \leqslant p, q < \infty$.
A Banach lattice $X$ is said to be \emph {$p$-convex} with constant $C$ if
$
\left\| \left( \sum_{j = 1}^N |f_j|^p \right)^{\frac 1 p} \right\|_X \leqslant C \left( \sum_{j = 1}^n \|f_j\|_X^p \right)^{\frac 1 p}
$
for any $\{f_j\}_{j = 1}^N \subset X$;
lattice $X$ is said to be \emph {$q$-concave} with constant $c$ if
$
\left( \sum_{j = 1}^n \|f_j\|_X^q \right)^{\frac 1 q} \leqslant c \left\| \left( \sum_{j = 1}^N |f_j|^q \right)^{\frac 1 q} \right\|_X
$
for any $\{f_j\}_{j = 1}^N \subset X$.
If $X$ is $p$-convex then $X'$ is $p'$-concave, and if $X$ is $q$-concave then $X'$ is $q'$-convex.
It is well known (see, e.~g., \cite [Book~II, Proposition~1.d.8] {cbs})
that a Banach lattice that is $p$-convex and $q$-concave can be renormed to make both its $p$-convexity and $q$-concavity constants
$C = c = 1$.
\begin {abridged}
It is easy to see that a lattice with the Fatou property which is $p$-convex and $q$-concave with some $1 < p, q < \infty$ is reflexive and
has order continuous norm.
\end {abridged}
\begin {unabridged}
The assumption of $p$-convexity imposed on a lattice $X$ enables us to raise $X$ to a power $1 < p < \infty$
without it becoming quasi-Banach
since $p$-convexity of $X$ is equivalent to $1$-convexity of $Y = X^p$.
This in turn implies that $X = Y^{\frac 1 p}$ and $X' = (Y')^{\frac 1 p} \lclassg {1}^{1 - \frac 1 p}$ provided that $X$ has the Fatou property,
so in this case $X'$ has order continuous norm and therefore $X = X'' = X'^*$.  By the same argument, if a lattice $X$ has the Fatou property and
$X$ is $q$-concave for some $1 < q < \infty$ then $X$ has order continuous norm and $X' = X^*$.
Thus a lattice $X$ which is both $p$-convex and $q$-concave with some $1 < p, q < \infty$
is reflexive, and also both $X$ and $X'$ have order continuous norm and enjoy many other nice properties.
\end {unabridged}

For a quasi-normed lattice $X$ and weights $\weightw$ such that $0 \leqslant \weightw \leqslant \infty$ almost everywhere
the weighted lattice $X (\weightw)$
is defined by
$$
X (\weightw) = 
\left\{g \mid \frac g {\weightw} \in X\right\}
$$
with the quasi-seminorm defined by
$\|f\|_{X (\weightw)} = \|f \weightw^{-1}\|_X$.  This somewhat cumbersome definition is needed because
the more natural definition
$
X (\weightw) = \{\weightw h \mid h \in X\}
$
is meaningless if the weight $\weightw$ takes value $+\infty$ on a set of positive measure
and it seems to be easier to allow this in the definition and work with weighted lattices that may be quasi-normed rather than negotiate
finiteness of $\weightw$ every time.
Thus in this setting one has $g = 0$ on the set where $\weightw = 0$, $g$ restricted on the set $\{\weightw = +\infty\}$
is an arbitrary measurable function,
and $\|\cdot\|_{X (\weightw)}$ is a norm for weights $\weightw$
such that $(\nu \times \mu) \left(\{\weightw = +\infty\}\right) = 0$.
If $\weightw = 0$ on a set of positive measure, we regard
$X (\weightw)$ as merely a set of functions with a seminorm under our conventions,
since then $\supp X (\weightw) \neq \supp X$.
In majorization arguments it is usually possible to avoid dealing with ``bad'' weights with the help of the following proposition.
\begin {proposition} [{\cite [Proposition~3.2] {rutsky2011en}}]
\label {xpos}
Suppose that $X$ is a Banach lattice on $(\Sigma, \mu)$.
Then for every $f \in X$ such that $f \neq 0$ identically and $\varepsilon > 0$ there exists $g \in X$
such that $g > |f|$ a. e. and $\|g\|_X \leqslant (1 + \varepsilon) \|f\|_X$.
\end {proposition}
The construction of a weighted lattice
yields
$$
\lclass {\infty} {\weightw} = \{ f \mid |f| \leqslant C \weightw \text { a. e. } \}.
$$
It is easy to see that $[X (\weightw)]' = X' (\weightw^{-1})$.
Notice that this definition of the weighted Lebesgue space $\lclass {p} {\weightw}$ differs from the ``classical'' one
with the norm defined by
$\|f\|_{p, \weightw}^p = \int |f|^p \weightw$, which is often used in the literature; the latter norm corresponds to the norm of the lattice
$\lclass {p} {\weightw^{-\frac 1 p}}$
in our notation.
Thus all weighted lattices are defined in the same way everywhere in this paper; however, one has to pay attention to this difference.
We adopt the natural conventions
$0^{-1} = \infty$ and $\infty^{-1} = 0$ in all expressions involving weights.

\begin {abridged}
The (centered) Hardy-Littlewood maximal operator
$$
M f (x, t) 
= \sup_{r > 0} \frac 1 {\nu (B (x, r))} \int_{B (x, r)} |f (z, t)| d\nu (z),
\quad x \in S, \quad t \in \Omega,
$$
is well-defined for a. e. $x \in S, t \in \Omega$, and the measurable functions $f$ on $(S \times \Omega, \nu \times \mu)$
that are locally summable in the first variable.
We say that a non-negative measurable function $\weightw$ on $(S \times \Omega, \nu \times \mu)$
belongs to the Muckenhoupt class $\apclass {p}$ for some $1 \leqslant p < \infty$
with a constant $C$ if
$$
\esssup_{t \in \Omega} \|M\|_{\lclass {p} {\weightw^{-1/p} (\cdot, t)} \to \lclass {p, \infty} {\weightw^{-1/p} (\cdot, t)}} \leqslant C.
$$
In the case $p > 1$ this condition is equivalent to
$$
\esssup_{t \in \Omega} \|M\|_{\lclass {p} {\weightw^{-1/p} (\cdot, t)}} \leqslant C'
$$
with a constant $C'$ estimated in terms of $C$ and $p$.
The class $A_1$ is characterized by the estimate $M \weightw \leqslant C' \weightw$ almost everywhere,
while classes $\apclass {p}$ for $p > 1$ are characterized by the well-known Muckenhoupt condition
\begin{multline}
\label {muckenhouptc}
\esssup_{x \in S, t \in \Omega} \,\, \sup_{r > 0} \bigg[\frac 1 {\nu 
\left(B (x, r)\right)} \int\limits_{B (x, r)} \weightw (u, t) d\nu (u) 
\bigg] 
\\
\quad\times\bigg[\frac 1 {\nu \left(B (x, r)\right)} \int\limits_{B (x, r)} \weightw (u, t)^{-\frac 1 {p - 1}} d\nu (u) \bigg]^{p - 1} 
< \infty.
\end{multline}

\begin {definition}
\label {aprdef}
A quasi-normed lattice $X$ on $(S \times \Omega, \nu \times \mu)$ is $\apclass {p}$-regular with constants $(C, m)$ if
for any $f \in X$ there exists a majorant $g \in X$, $g \geqslant |f|$ such that $\|g\|_X \leqslant m \|f\|_X$
and $g \in \apclass {p}$ with constant $C$.
\end {definition}

\begin {unabridged}
The following simple proposition justifies the title of this paper.
\end {unabridged}
\begin {proposition} [{\cite [Proposition~1.2] {rutsky2011en}}]
\label {mbnecc}
A quasi-normed lattice $X$ on $(S \times \Omega, \nu \times \mu)$ is $\apclass {1}$-regular if and only if
the maximal operator $M$ is bounded in $X$.
\end {proposition}

It is easy to see that the $\apclass {1}$-regularity property is self-improving,
which is the subject of the following proposition; see also \cite {lernerperez2007}.
\begin {proposition}
\label {a1si}
Suppose that $X$ is an $\apclass {1}$-regular Banach lattice on $(S \times \Omega, \nu \times \mu)$ with constants $(C, m)$.
Then $X^r$ is also an $\apclass {1}$-regular lattice for some $r > 1$ depending only on $C$.
\end {proposition}
Indeed, let $r > 1$ be the constant of the reverse H\"older inequality that is satisfied for all $\apclass {1}$ weights with constant $C$.
Suppose that $f \in X^r$, and let $g$ be an $\apclass {1}$-majorant for $|f|^{\frac 1 r}$ in $X$ with constants $(C, m)$.
Then $g^r$ is an $\apclass {1}$-majorant for $f$ with constants independent of $f$, because by the reverse H\"older inequality we have
an estimate
\begin {multline*}
\frac 1 {\nu (B (x, \rho))} \int_{B (x, \rho)} g^r (u, \omega) d\nu (u) \leqslant
\\
c \left( \frac 1 {\nu (B (x, \rho))} \int_{B (x, \rho)} g (u, \omega) d\nu (u) \right)^r \leqslant
c\, C^r [g (x, \omega)]^r
\end {multline*}
for almost all $x \in S$, $\omega \in \Omega$ and $\rho > 0$ with a constant $c$ independent of $f$, $x$, $\omega$ and $\rho$.

We say that $T$ is a \emph {Calderon-Zygmund operator} if $T$ is a singular integral operator that is bounded in $\lclass {2} {\mathbb R^n}$,
its kernel $K (x, y)$ satisfies
\begin {equation}
\label {kernelest}
|K (x, s) - K (x, t)| \leqslant C_K \frac {|s - t|^\gamma} {|x - s|^{n + \gamma}}, \quad x, s, t \in \mathbb R^n, \quad |x - s| > 2 |s - t|
\end {equation}
with some $\gamma > 0$,
and the kernel $K^* (y, x) = K (x, y)$ of the adjoint operator $T^*$ satisfies the same estimates.  It is well known that $T$ is bounded in $\lclassg {p}$
for all $1 < p < \infty$.

The Coifman-Fefferman inequality \cite {coifmanfefferman1974}
\begin {equation}
\label {coifmani}
\int |T f|^p \omega \leqslant C \int (M f)^p \omega, \quad 0 < p < \infty,
\end {equation}
holds true for Calderon-Zygmund operators $T$ with $C$ independent of $f$
and any weight $\omega \in \apclass {\infty}$
with all locally summable functions $f$ such that the right-hand part of \eqref {coifmani} is finite.

\begin {proposition}
\label {czobounda1p}
Suppose that $X$ is a Banach lattice on $\mathbb R^n \times \Omega$ having either the Fatou property
or order continuous norm
and both $X$ and $X'$ are $\apclass {1}$-regular.
Then any Calderon-Zygmund operator $T$ is bounded in~$X$.
\end {proposition}
Indeed, let $f \in X$ and $g \in X'$, and let $h$ be an $\apclass {1}$-majorant of $g$ in $X'$.
Then
$$
\int (M f) h \leqslant \|M f\|_X \|h\|_{X'} \leqslant c_1  \|f\|_X \|g\|_{X'} < \infty,
$$
and
the Coifman-Fefferman inequality \eqref {coifmani} with $p = 1$ implies that
\begin {equation*}
\int (T f) g \leqslant \int |T f| h \leqslant c \int (M f) h \leqslant c\, c_1  \|f\|_X \|g\|_{X'}
\end {equation*}
with certain constants $c$ and $c_1$ independent of $f$ and $g$, which shows that $T$ acts boundedly in $X$.

\begin {definition}
\label {nondego}
A mapping $T : \lclassg {2} \to \lclassg {2}$ is called $\apclass {2}$-nondegenerate with a constant $C$ if
boundedness of $T$ in a lattice $\lclass {2} {\weightw^{-\frac 1 2}}$ implies $\weightw \in \apclass {2}$ with constant $C$.
\end {definition}

We remark that by \cite [Proposition~3.7] {rutsky2011en}
an $\apclass {2}$-nondegenerate linear operator $T$ in $\lclass {2} {\mathbb R^n}$
is also $\apclass {2}$-nondegenerate as an operator in $\lclass {2} {\mathbb R^n \times \Omega}$ acting
in the first variable.

\begin {proposition} [{\cite [Chapter~5, \S4.6] {stein1993}}]
\label {tndcs}
Suppose that $T$ is a Calderon-Zygmund operator with kernel $K$
and there exist some $u \in \mathbb R^n$ and a constant $c$ such that for any $x \in \mathbb R^n$ and $t \neq 0$ we have
\begin {equation}
\label {sdir}
|K (x, x + t u)| \geqslant c t^{-n}.
\end {equation}
Then $T$ is $\apclass {2}$-nondegenerate.
\end {proposition}
It is easy to see that the Hilbert transform $H$ on $\mathbb R$ with kernel $K (x, y) = \frac {c_1} {x - y}$
and Riesz transforms $R_j$, $1 \leqslant j \leqslant n$ on $\mathbb R^n$
with kernels $K_j (x, y) = \frac {c_n (y_j - x_j)} {|y - x|^{n + 1}}$, where $c_n \neq 0$ are some constants,
satisfy condition~\eqref {sdir} for $u = e_j$, $e_j$ being the $j$-th coordinate basis vector of $\mathbb R^n$,
and thus all these operators are $\apclass {2}$-nondegenerate.
\end {abridged}
%\end {abridged}

%\begin {unabridged}
\begin {comment}
In this section we introduce some basic notions.
Suppose that $(S, \nu)$ is a space of homogeneous type, i.~e. $S$ is
a quasimetric space equipped with a Borel measure $\nu$ that has the doubling property:
$$
\nu (B (x, 2 r)) \leqslant c \nu (B (x, r))
$$
for all $x \in S$ and $0 < r < \infty$ where $c$ is a constant called a
doubling constant.  Futher information on these spaces and real harmonic analysis on them see e. g. in  \cite {denghan, stein1993}.
It is customary in the literature on real harmonic analysis to consider only the typical case of $S$ being a Euclidean space $\mathbb R^n$
equipped with the Lebesgue measure.
However, for the main results of this paper this restriction is not at all necessary.
We are working with (real or complex) quasi-normed lattices of measurable functions on
$(S \times \Omega, \nu \times \mu)$
where $(\Omega, \mu)$ is an arbitrary measurable space with a $\sigma$-finite measure $\mu$.
The second variable $t \in \Omega$ allows us to treat vector-valued results in this setting
and it does not usually cause much trouble in the arguments.
The case of just one variable $x$ is, of course, included in this generality: one just takes the point mass measure for $\mu$.
A quasi-normed lattice of measurable functions $X$ is by definition a quasi-normed space of measurable functions $X$
in which the norm is compatible with the natural order; that is, if $|f| \leqslant g$ a.~e. for some function $g \in X$
then $f \in X$ and $\|f\|_X \leqslant C \|g\|_X$ for some constant $C$ independendent of $f$ and $g$. 
One usually has $C = 1$; for example, it is always the case when the lattice has the Fatou property which will be introduced shortly.
Further information on lattices can be found in e. g. \cite [Chapter~10]{kantorovichold}.
For simplicity we only work with lattices $X$ such that $\supp X = S \times \Omega$.
Many interesting normed spaces appearing in analysis are or can be represented as such lattices; for example
the Lebesgue spaces $\lclassg {p}$, the Orlicz spaces $\lclassg {M}$, and 
the Lebesgue spaces $\lclassg {p (\cdot)}$ with variable exponent $p (\cdot)$ along with general modular spaces.
The second variable allows us to naturally include lattices with mixed norm such as $\lclass {p} {\lclassg {r}}$ in this setting.

Now we introduce certain properties of quasi-normed lattices and some related objects.
For a Banach lattice of measurable functions $X$, any order continuous functional $f$ on $X$ (order continuity means that given a sequence $x_n \in X$
such that $\sup_n |x_n| \in X$ and $x_n \to 0$ a. e. one also has $f (x_n) \to 0$) has an integral representation
$f (x) = \int x y_f$ for some measurable function $y_f$ which can be identified with $f$.
The set of all such functionals $X'$ is a Banach lattice with
the norm defined by $\|f\|_{X'} = \sup_{g \in X, \|g\|_X = 1} \int |f g|$.  The lattice $X'$ is called the order dual of the lattice $X$.
The norm of a lattice $X$ is said to be order continuous if for any nonincreasing sequence
$x_n \in X$ converging to 0 a. e. one also has
$\|x_n\|_X \to 0$.  The norm of a Banach lattice $X$ is order continuous if and only if $X^* = X'$, i.~e. if an only if all norm continous functionals on
$X$ are order continous.
A lattice $X$ has the Fatou property if for any
$f_n, f \in X$ such that $\|f_n\|_X \leqslant 1$ and the sequence $f_n$ 
converges to $f$ a. e. it is also true that $f \in X$ and $\|f\|_X \leqslant 1$.
The Fatou property of a lattice $X$ is equivalent to $(\nu \times \mu)$-closedness of the unit ball
$B_X$ of the lattice $X$
(here and elsewhere $(\nu \times \mu)$-convergence denotes the
convergence in measure in any measurable set $E$ such that $(\nu \times \mu) (E) < \infty$).
If the lattice $X$ is Banach then the Fatou property is equivalent to order reflexivity of $X$, i.~e. to the relation $X'' = X$.
For a lattice $X$ either one of the Fatou property or the order continuity of norm property %of a lattice $X$
is sufficient to guarantee that the lattice $X'$ is a norming set of functionals for
$X$, i.~e. that $\|f\|_X = \sup_{g \in X', \|g\|_{X'} = 1} \int f g$ for all $f \in X$. The order dual $X'$ of a Banach lattice $X$ always has the Fatou property.

We now introduce a couple of very useful constructions that play a major part in the analysis of lattices and their majorization
and other properties.
For any two quasi-normed lattices $X$ and $Y$ on the same measurable space
the set of poinwise products of their functions
$$
X Y = \{ f g \mid f \in X,\break g \in Y\}
$$
is a quasi-normed lattice with the norm defined by
$$
\|h\|_{X Y} = \inf_{h = f g} \|f\|_X \|g\|_Y.
$$
If both lattices $X$ and $Y$ satisfy the Fatou property then the lattice $X Y$ also has the Fatou property.
This lattice multiplication is associative: if $X$, $Y$ and $Z$ are lattices of measurable functions on the same measurable space then
$(X Y) Z = X (Y Z)$ (here and elsewhere, if not stated otherwise,
the equality of lattices is understood as the equality of the sets together with the equality of the quasi-norms).
It is easy to see that if either of the lattices $X$ and $Y$ has order continuous quasi-norm then the norm of the lattice $X Y$ is also order continuous.

For any $\delta > 0$ and a quasi-normed lattice $X$ the lattice $X^\delta$ consists of all measurable functions $f$
having well-defined and finite quasi-norms $\|f\|_{X^\delta} = \| |f|^{1/\delta} \|_X^\delta$.
For example, $\lclassg {p}^\delta = \lclassg {\frac p \delta}$.
For any product $X Y$ of quasi-lattices and $\delta > 0$ one has the natural relation $(X Y)^\delta = X^\delta Y^\delta$.
If a lattice $X$ has the Fatou property then $X^\delta$ also has the Fatou property;
if a lattice $X$ has order continous quasi-norm then $X^\delta$ also has order continuous quasi-norm.
For any $0 < \delta \leqslant 1$ if $X$ is a Banach lattice then
$X^\delta$ is also a Banach lattice.
If $X$ and $Y$ are Banach lattices then for any $0 < \delta < 1$ lattice
$X^{1 - \delta} Y^\delta$, sometimes called the \emph {Calderon-Lozanovsky product} of $X$ and $Y$, is also Banach; moreover, one has a very useful relation
$(X^{1 - \delta} Y^\delta)' = (X')^{1 - \delta} (Y')^\delta$ (see \cite {calderon1964}, \cite {lozanovsky1969}).
If $Z = X^{1 - \delta} Y^\delta$ has either the Fatou property or has order continuous norm
then 
$Z$ is an exact interpolation space of exponent $\delta$ between $X$ and $Y$
(which coincides with either the complex interpolation space $(X, Y)^\theta$ or with $(X, Y)_\theta$ respectively), i.~e. any linear operator
$T$ that acts boundedly in $X$ and in $Y$ also acts boundedly in $Z$ and $\|T\|_Z \leqslant \|T\|_X^{1 - \delta} \|T\|_Y^\delta$;
see, e.~g., \cite {lozanovsky1972}, \cite {calderon1964}, \cite {kps} for more detail.

Let $1 \leqslant p, q < \infty$.
A Banach lattice $X$ is said to be \emph {$p$-convex} with constant $C$ if
$$
\left\| \left( \sum_{j = 1}^N |f_j|^p \right)^{\frac 1 p} \right\|_X \leqslant C \left( \sum_{j = 1}^n \|f_j\|_X^p \right)^{\frac 1 p}
$$
for any $\{f_j\}_{j = 1}^N \subset X$;
$X$ is said to be \emph {$q$-concave} with constant $c$ if
$$
\left( \sum_{j = 1}^n \|f_j\|_X^q \right)^{\frac 1 q} \leqslant c \left\| \left( \sum_{j = 1}^N |f_j|^q \right)^{\frac 1 q} \right\|_X
$$
for any $\{f_j\}_{j = 1}^N \subset X$.
If $X$ is $p$-convex then $X'$ is $p'$-concave, and if $X$ is $q$-concave then $X'$ is $q'$-convex.
It is well known (see, e.~g., \cite [Book~II, Proposition~1.d.8] {cbs})
that a Banach lattice that is $p$-convex and $q$-concave can be renormed so that its $p$-convexity and $q$-concavity constants
are both $1$.  The assumption of $p$-convexity imposed on a lattice $X$ enables us to raise the lattice in a power $1 < p < \infty$
without it becoming quasi-Banach
since $p$-convexity of $X$ is equivalent to $1$-convexity of $Y = X^p$.
This in turn implies that $X = Y^{\frac 1 p}$ and $X' = (Y')^{\frac 1 p} \lclassg {1}^{1 - \frac 1 p}$ provided that $X$ has the Fatou property,
so in this case $X'$ has order continuous norm and therefore $X = X'' = X'^*$.  By the same argument if a lattice $X$ has the Fatou property and
$X$ is $q$-concave for some $1 < q < \infty$ then $X$ has order continuous norm and $X' = X^*$.
Thus a lattice $X$ which is both $p$-convex and $q$-concave with some $1 < p, q < \infty$
is reflexive, and also both $X$ and $X'$ have order continuous norm and enjoy many other nice properties.

%For a lattice of measurable functions $X$ by $X_+$ we denote the cone of its positive functions $\{f \in X \mid f \geqslant 0 \text { a. e.} \}$.

For a quasi-normed lattice $X$ and weights $\weightw$ such that $0 \leqslant \weightw \leqslant \infty$ almost everywhere
the weighted lattice $X (\weightw)$
is defined by
$$
X (\weightw) = 
\left\{g \mid \frac g {\weightw} \in X\right\}
$$
with the quasi-seminorm defined by
$\|f\|_{X (\weightw)} = \|f \weightw^{-1}\|_X$.  This somewhat cumbersome definition is needed because
the more natural definition
$$
X (\weightw) = \{\weightw h \mid h \in X\}
$$
is meaningless if the weight $\weightw$ takes value $+\infty$ on a set of positive measure
and it seems to be easier to allow this in the definition and work with weighted lattices that may be quasi-normed rather than negotiate
finiteness of $\weightw$ every time.
Thus in this setting one has $g = 0$ on the set where $\weightw = 0$, $g$ restricted on the set $\{\weightw = +\infty\}$
is an arbitrary measurable function,
and $\|\cdot\|_{X (\weightw)}$ is a norm for weights $\weightw$
such that $\nu \times \mu \left(\{\weightw = +\infty\}\right) = 0$.
If $\weightw = 0$ on a set of positive measure, we regard
$X (\weightw)$ as merely a set of functions with a seminorm under our conventions, since then $\supp X (\weightw) \neq \supp X$.
Note that in majorization arguments it is usually possible to avoid dealing with ``bad'' weights with the help of the following proposition.
\begin {proposition} [{\cite [Proposition~3.2] {rutsky2011en}}]
\label {xpos}
Suppose that $X$ is a Banach lattice on $(\Sigma, \mu)$.
Then for every $f \in X$ such that $f \neq 0$ identically and $\varepsilon > 0$ there exists $g \in X$
such that $g > |f|$ a. e. and $\|g\|_X \leqslant (1 + \varepsilon) \|f\|_X$.
\end {proposition}
The construction of a weighted lattice
yields
$$
\lclass {\infty} {\weightw} = \{ f \mid |f| \leqslant C \weightw \text { a. e. } \}.
$$
It is easy to see that $[X (\weightw)]' = X' (\weightw^{-1})$.
Notice that this definition of the weighted Lebesgue space $\lclass {p} {\weightw}$ differs from the ``classical'' one
with the norm defined by
$\|f\|_{p, \weightw}^p = \int |f|^p \weightw$, which is often used in the literature; the latter norm corresponds to the norm of the lattice
$\lclass {p} {\weightw^{-\frac 1 p}}$
in our notation.
Thus all weighted lattices are defined in the same way everywhere in this paper; however, one has to pay attention to this difference.
We adopt the natural conventions
$0^{-1} = \infty$ and $\infty^{-1} = 0$ in all expressions involving weights.
%\end {unabridged}
\end {comment}

\begin {unabridged}

\section {Muckenhoupt weights and $\apclass {p}$-majorants}

\label {mwaam}

In this section we introduce some useful notions having to do with the Muckenhoupt weights;
for more detail see, e.~g., \cite [Chapter~5] {stein1993}.
The (centered) Hardy-Littlewood maximal operator
$$
M f (x, t) 
= \sup_{r > 0} \frac 1 {\nu (B (x, r))} \int_{B (x, r)} |f (z, t)| d\nu (z),
\quad x \in S, \quad t \in \Omega,
$$
is well-defined for a. e. $x \in S, t \in \Omega$, and the measurable functions $f$ on $(S \times \Omega, \nu \times \mu)$
that are locally summable in the first variable.
We say that a non-negative measurable function $\weightw$ on $(S \times \Omega, \nu \times \mu)$
belongs to the Muckenhoupt class $\apclass {p}$ for some $1 \leqslant p < \infty$
with a constant $C$ if
$$
\esssup_{t \in \Omega} \|M\|_{\lclass {p} {\weightw^{-1/p} (\cdot, t)} \to \lclass {p, \infty} {\weightw^{-1/p} (\cdot, t)}} \leqslant C.
$$
In the case $p > 1$ this condition is equivalent to
$$
\esssup_{t \in \Omega} \|M\|_{\lclass {p} {\weightw^{-1/p} (\cdot, t)}} \leqslant C'
$$
with a constant $C'$ estimated in terms of $C$ and $p$.
The class $A_1$ is characterized by the estimate $M \weightw \leqslant C' \weightw$ almost everywhere,
while classes $\apclass {p}$ for $p > 1$ are characterized by the well-known Muckenhoupt condition
\begin{multline}
\label {muckenhouptc}
\esssup_{x \in S, t \in \Omega} \,\, \sup_{r > 0} \bigg[\frac 1 {\nu 
\left(B (x, r)\right)} \int\limits_{B (x, r)} \weightw (u, t) d\nu (u) 
\bigg] 
\\
\quad\times\bigg[\frac 1 {\nu \left(B (x, r)\right)} \int\limits_{B (x, r)} \weightw (u, t)^{-\frac 1 {p - 1}} d\nu (u) \bigg]^{p - 1} 
< \infty.
\end{multline}
 %see e. g. \cite [Chapter~5] {stein1993}.
The class $\apclass {\infty}$ is defined as the class of weights $\weightw$ satisfying the reverse H\"older inequality
\begin{multline}
\label {rhi}
\esssup_{x \in S, t \in \Omega} \,\, \sup_{r > 0} \bigg[\frac 1 {\nu 
\left(B (x, r)\right)} \int\limits_{B (x, r)} [\weightw (u, t)]^q d\nu (u) 
\bigg]^{\frac 1 q}
\\
\quad\times\bigg[\frac 1 {\nu \left(B (x, r)\right)} \int\limits_{B (x, r)} \weightw (u, t) d\nu (u) \bigg]^{-1} 
< \infty
\end{multline}
with some $q > 1$, and for certainty we will take for the value of the supremum in \eqref {rhi}
for the $\apclass {\infty}$ constant of the weight $\weightw$.  It is well known that $\weightw \in \apclass {\infty}$ if and only if
$\weightw \in \apclass {p}$ with some $1 < p < \infty$ and the $\apclass {p}$ constant of the weight $\weightw$
depending only on the $\apclass {\infty}$ constant of the weight $\weightw$ and vice versa.

The following notion is a natural refinement of the $\BMO$-regularity property
which was apparently first introduced by N.~Kalton in \cite {kalton1994}.
\begin {definition}
\label {aprdef}
A quasi-normed lattice $X$ on $(S \times \Omega, \nu \times \mu)$ is $\apclass {p}$-regular with constants $(C, m)$ if
for any $f \in X$ there exists a majorant $g \in X$, $g \geqslant |f|$ such that $\|g\|_X \leqslant m \|f\|_X$
and $g \in \apclass {p}$ with constant $C$.
\end {definition}
This property was formally introduced and studied to some extent in \cite {rutsky2011en}; we will reference here only the results used in the present work.
 %The following simple proposition justifies the use of the term ``$\apclass {1}$-regularity''.
\begin {proposition} [{\cite [Proposition~1.2] {rutsky2011en}}]
\label {mbnecc}
A quasi-normed lattice $X$ on $(S \times \Omega, \nu \times \mu)$ is $\apclass {1}$-regular if and only if
the maximal operator $M$ is bounded in $X$.
\end {proposition}
Sufficiency is trivial, and necessity quickly follows from an application of the famous Rubio de Francia construction.
\begin {comment}
The following proposition is a direct consequence of duality and
the connection between Muckenhoupt weights and boundedness of the maximal operator.
\begin {proposition} [{\cite [Proposition~2.3] {rutsky2011en}}]
\label {aptoconj}
Suppose that $X$ is a Banach lattice on $(S \times \Omega, \nu \times \mu)$
such that $X'$ is a norming space for $X$.
If $X'$ is $\apclass {1}$-regular then $X^{\frac 1 q}$ is $\apclass {1}$-regular for all $q > 1$.
If $X'$ is $\apclass {p}$-regular with some $p > 1$ then
$X^{\frac 1 p}$ is $\apclass {1}$-regular.
\end {proposition}

 %The following result was obtained in \cite {rutsky2011en} in order to establish the division and self-duality property of $\BMO$-regularity
%in the general setting.
\begin {theorem} [{\cite [Theorem~1.6] {rutsky2011en}}]
\label {bmomdiv}
Suppose that $X$ is a Banach lattice on $(S \times \Omega, \nu \times \mu)$
having the Fatou property.  Suppose also that
$X \lclassg {q}$ with some $1 < q < \infty$ is a Banach lattice and $X \lclassg {q}$ is $\apclass {p}$-regular for some $1 \leqslant p < \infty$.
Then $X$ is $\apclass {p + 1}$-regular.
\end {theorem}

This result is a direct precursor to a very deep and nontrivial fact that the BMO-regularity property is self-dual;
 % at least
%for Banach lattices having the Fatou property;
see \cite {kalton1994}, \cite {kisliakov2002en}, \cite {rutsky2011en}.
The only currently known method of establishing Theorem~\ref {bmomdiv} (in the case of $\mathbb R^n$) is through an application of a fixed point theorem
somewhat similar to what is done in Section~\ref {aprmainlemma} below.
\end {comment}

As a consequence of the reverse H\"older inequality we see that the $\apclass {1}$-regularity property is self-improving,
which is the subject of the following proposition.
(It is not difficult to see that the general $\apclass {p}$-regularity property is also self-improving in this manner,
but we will not need it in the present work).
There is also a fairly general approach that makes it possible to establish this property using certain methods
originating in the geometry of Banach spaces;
see \cite {lernerperez2007}.
\begin {proposition}
\label {a1si}
Suppose that $X$ is an $\apclass {1}$-regular Banach lattice on $(S \times \Omega, \nu \times \mu)$ with constants $(C, m)$.
Then $X^r$ is also an $\apclass {1}$-regular lattice for some $r > 1$ depending only on $C$.
\end {proposition}
Indeed, let $r > 1$ be a constant of the reverse H\"older inequality that is satisfied for all $\apclass {1}$ weights with constant $C$.
Suppose that $f \in X^r$, and let $g$ be an $\apclass {1}$-majorant for $|f|^{\frac 1 r}$ in $X$ with constants $(C, m)$.
Then $g^r$ is an $\apclass {1}$-majorant for $f$ with constants independent of $f$, because by the reverse H\"older inequality we have
an estimate
\begin {multline*}
\frac 1 {\nu (B (x, \rho))} \int_{B (x, \rho)} g^r (u, \omega) d\nu (u) \leqslant
\\
c \left( \frac 1 {\nu (B (x, \rho))} \int_{B (x, \rho)} g (u, \omega) d\nu (u) \right)^r \leqslant
c\, C^r [g (x, \omega)]^r
\end {multline*}
for almost all $x \in S$, $\omega \in \Omega$ and $\rho > 0$ with a constant $c$ independent of $f$, $x$, $\omega$ and $\rho$.

\section {Calderon-Zygmund operators}

\label {czos}

In this section we will show how certain conditions on the lattices are
sufficient for boundedness of the Calderon-Zygmund operators in the general setting.
Namely, we give 4 somewhat independent proofs of the implication $1 \Rightarrow 2$ of Theorem~\ref {themcr}
that use different properties of a Calderon-Zygmund operator.
For generalities on the real harmonic analysis see, e.~g., \cite {denghan}, \cite {stein1993}.
 %For a standard survey on the singular integral operators see, e.~g., \cite {stein1993}.
Although it seems possible to extend all of the results used here to the general setting of spaces of homogeneous type and beyond,
for simplicity we will only discuss the standard setting of $\mathbb R^n$ with the usual Lebesgue measure $d\nu = dm$.
We will also use, in contrast to the definition introduced in Section~\ref {mwaam}, the (uncentered) Hardy-Littlewood maximal operator defined by
$$
M f (x) = \sup_{Q \ni x} \frac 1 {|Q|} \int_Q |f (y)| dy, \quad x \in \mathbb R^n,
$$
where the supremum is taken over all cubes $Q \subset \mathbb R^n$ with sides parallel to the coordinate axes.
However, it is well known that this definition is pointwise equivalent %up to multiplicative constants
to the one given before.

We say that $T$ is a \emph {Calderon-Zygmund operator} if $T$ is a singular integral operator that is bounded in $\lclass {2} {\mathbb R^n}$
and its kernel $K (x, y)$ satisfies
\begin {equation}
\label {kernelest}
|K (x, s) - K (x, t)| \leqslant C_K \frac {|s - t|^\gamma} {|x - s|^{n + \gamma}}, \quad x, s, t \in \mathbb R^n, \quad |x - s| > 2 |s - t|
\end {equation}
with some $\gamma > 0$,
and the kernel $K (y, x)$ of the adjoint operator $T^*$ satisfies the same estimates.  It is well known that $T$ is bounded in $\lclassg {p}$
for all $1 < p < \infty$.
 %Although recently interesting new tools were developed for treatment of such operators in a fairly general setting,
We begin %set off
with the more classical approach, which is also very simple.
The following proposition contains the implication $1 \Rightarrow 2$ of Theorem~\ref {themcr}
(see also Corollary~\ref {czobounda1pc} in Section~\ref {aprmainlemma} below).
 %Using these two ingredients, it is easy to establish a sufficient condition for boundedness of $T$ on $X$
%that, as we will see later, is also often necessary. % FIXME: reference
\begin {proposition}
\label {czobounda1p}
Suppose that $X$ is a Banach lattice on $\mathbb R^n$ having either the Fatou property
or order continuous norm,
$X$ is $\apclass {1}$-regular and $X'$ is $\apclass {\infty}$-regular.
Then any Calderon-Zygmund operator $T$ is bounded in~$X$.
\end {proposition}
Indeed, let $f \in X$ and $g \in X'$, and let $h$ be an $\apclass {\infty}$-majorant of $g$ in $X'$.
Then
$$
\int (M f) h \leqslant \|M f\|_X \|h\|_{X'} \leqslant c_1  \|f\|_X \|g\|_{X'} < \infty,
$$
and
the Coifman-Fefferman inequality \eqref {coifmani} with $p = 1$ implies that
\begin {equation*}
\int (T f) g \leqslant \int |T f| h \leqslant c \int (M f) h \leqslant c\, c_1  \|f\|_X \|g\|_{X'}
\end {equation*}
with certain constants $c$ and $c_1$ independent of $f$ and $g$, which implies that $T$ acts boundedly in $X$.
Compared to the other approaches that follow (that, at least in their presently available form,
significantly rely in their details on the structure of the dyadic cubes in $\mathbb R^n$
which makes it harder to carry the arguments over to a general space of homogeneous type),
it is easy to see that the Coifman-Fefferman inequality and
other parts of the proof remain valid in the general case of Calderon-Zygmund operators
on $\sigma$-finite spaces of measurable functions on $S \times \Omega$ where $S$ is a space of homogeneous type.

Now let us briefly describe another approach to establishing %the implication $1 \Rightarrow 2$ of Theorem~\ref {themcr}
a slightly weaker version of Proopsition~\ref {czobounda1p}
that uses the classical Fefferman-Stein inequality that was recently generalized to general Banach lattices.
The Fefferman-Stein maximal function $f^\sharp$ on $\mathbb R^n$ is defined for a locally integrable function $f$ by
$$
f^\sharp (x) = \sup_{Q \ni x} \frac 1 {|Q|} \int_Q \left|f (y) - f_Q \right| dy, \quad x \in \mathbb R^n,
$$
where the supremum is taken over all cubes $Q \subset \mathbb R^n$ 
containing $x$ with sides parallel to the coordinate axes
and $f_Q = \frac 1 {|Q|} \int_Q f (z) dz$ is the average of $f$ over $Q$ with respect to the Lebesgue measure.
This maximal function is very useful in the estimates of the Calderon-Zygmund operators $T$.
On the one hand, we have the well-known (and rather simple) pointwise estimate
\begin {equation}
\label {fsharpczo}
(T f)^\sharp \leqslant c_r \left( M |f|^r \right)^{\frac 1 r}
\end {equation}
almost everywhere for any $1 < r < \infty$;
see, e.~g., \cite [Chapter~4, \S4.2] {stein1993}.
On the other hand, there is the classical Fefferman-Stein inequality
$$
\|f\|_{\lclassg {p}} \leqslant c \|f^\sharp\|_{\lclassg {p}}
$$
for $1 < p < \infty$.
The latter was recently generalized to general Banach lattices as follows
(see \cite {lerner2010} for more information).
In the rather convenient notation
$S_0$ denotes the set of all measurable functions $f$ on $\mathbb R^n$ such that their nonincreasing rearrangement
$f^*$ satisfies $f^* (+\infty) = 0$, and the main tools of \cite {lerner2010}
that will appear shortly in this seciton work for this class of functions rather than just
the locally summable ones.
Surely $S_0$ contains all measurable functions with compact support, and
thus it is easy to see that $S_0$ is dense in a Banach lattice $X$ if, for example,
$X$ has order continuous norm.  The converse, however, is not true:
take, for example, $X = \lclass {\infty} {\weightw}$ with a weight $\weightw$ satisfying $\weightw^* (+\infty) = 0$.
\begin {theorem} [{\cite [Corollary~4.3] {lerner2010}}]
\label {fsharpdual}
Suppose that $X$ is an $\apclass {1}$-regular real Banach lattice of measurable functions on $\mathbb R^n$
having the Fatou property.  Then the following conditions are equivalent.
\begin {enumerate}
\item
$X'$ is $\apclass {1}$-regular.
\item
There exists some $c > 0$ such that $\|f\|_X \leqslant c \|f^\sharp\|_X$ for all $f \in S_0 \cap X$.
\end {enumerate}
\end {theorem}
These two ingredients allow us to easily establish Proposition~\ref {czobounda1p} under the additional assumption that
$S_0$ is dense in $X$.
Indeed, by the assumed density property it is sufficient to estimate $\|T f\|_X$
for all $f \in S_0 \cap X$.  An application of Theorem~\ref {fsharpdual}, \eqref {fsharpczo} and Proposition~\ref {a1si}
yields
\begin {multline*}
\|T f\|_X \leqslant c \left\|(T f)^\sharp\right\|_X \leqslant c\, c_r \left\|\left(M |f|^r\right)^{\frac 1 r} \right\|_X \leqslant
\\
c\, c_r \|M\|_{X^r \to X^r}^{\frac 1 r} \|f\|_X \leqslant c_1 \|f\|_X
\end {multline*}
with some $r > 1$ and constants $c$, $c_1$ and $c_r$ independent of $f$.  %The proof of Proposition~\ref {czobounda1p} is complete.
The assumption in Theorem~\ref {fsharpdual} that $X$ is a real Banach lattice is easy to lift; see, e.~g.,
\cite [Proposition~6] {rutsky2013iv1}.

We need some more preliminaries before further discussion.
The \emph {Str\"omberg local sharp maximal function} is defined for $f \in S_0$ by
\begin {equation}
\label {sharpfdef}
M^\sharp_\lambda f (x) = \sup_{Q \ni x} \inf_{c \in \mathbb R} ((f - c)\chi_Q)^* (\lambda |Q|), \quad x \in \mathbb R^n,
\end {equation}
where the supremum is taken over all cubes $Q \subset \mathbb R^n$ containing $x$.
Functions $M^\sharp_\lambda f$ and $f^\sharp$ are closely related via the following estimate
which holds true
with all sufficiently small $0 < \lambda < 1$ and some $c_0, c_1 > 1$
for all locally summable $f$:
\begin {equation}
\label {sharpequiv}
c_0 M M^\sharp_\lambda f (x) \leqslant f^\sharp (x) \leqslant c_1 M M^\sharp_\lambda f (x), \quad x \in \mathbb R^n;
\end {equation}
see, e.~g., \cite {jawerthtorchinsky1985}, \cite {lerner2003}.
On the other hand, $M^\sharp_\lambda$ in many cases provides estimates that are significantly finer than those obtained by the means of the
Fefferman-Stein sharp maximal function.
For example (see \cite {alvarezperez1994}, \cite {jawerthtorchinsky1985}),
\begin {equation}
\label {flsharpczo}
M^\sharp_\lambda (T f) \leqslant c M f
\end {equation}
almost everywhere for all locally summable functions $f$ with $c$ independent of $f$.
Estimate \eqref {flsharpczo} is sharper than \eqref {fsharpczo}, as it corresponds to the missing
limiting case $r = 1$ in \eqref {fsharpczo}, and we can easily obtain \eqref {fsharpczo} from \eqref {flsharpczo}
using \eqref {sharpequiv} and \cite [Chapter~5, \S5.2] {stein1993}.
Finally, there is the following result similar to the well-known duality relation between $\hclassg {1}$ and $\BMO$.
\begin {theorem} [{\cite [Theorem~1] {lerner2004}}]
\label {mmsharpd}
$$
\int |f g| \leqslant c \int M^\sharp_\lambda f \, M g
$$
for any $f \in S_0$ and locally summable function $g$
with some $c$ and $\lambda$ independent of $f$ and $g$.
\end {theorem}

Examining the details of the previous argument
 %proof of Proposition~\ref {czobounda1p}
contained in Theorem~\ref {fsharpdual} quickly leads to
 %a more direct and essentially simplified proof of Proposition~\ref {czobounda1p} and even generalize it
the following observation.
 %to the general case of operators acting between different Banach lattices.
\begin {theorem}
\label {dpmxy}
Suppose that $X$, $Y$ and $Z$ are Banach lattices on $\mathbb R^n$ having the Fatou property,
 %the set $S_{00}$ of bounded functions with compact support is dense\footnote {Again,
 %if the norm of $X$ is order continuous then $S_{00}$ is dense in $X$
 %but not vice versa: take, for example, a variable exponent Lebesgue space $\lclassg {p (\cdot)}$ with $p = \infty$
 %on a certain bounded compact set $B$ and $1 \leqslant p < \infty$ elsewhere.}
 %in $X$,
$S_0$ is dense in $X$,
and suppose that the Hardy-Littlewood
maximal operator $M$ acts boundedly from $X$ to $Z$ and from $Y'$ to $Z'$.
Then any %Calderon-Zygmund
operator $T$ that satisfies estimate \eqref {flsharpczo} acts boundedly from $X$ to $Y$.
 % and the same is true for the
%maximal truncation operator $T_\natural$.
\end {theorem}
Theorem~\ref {dpmxy} 
follows at once from Theorem~\ref {mmsharpd}, since for any $f \in X \cap S_0$ and $g \in Y'$ we have
an estimate
\begin {multline}
\label {dpmxye}
\int |(T f) g| \leqslant c \int \left[M^\sharp_\lambda (T f)\right] M g \leqslant
c_1 \int (M f) (M g) \leqslant
\\
c_1 \| M f \|_Z \| M g \|_{Z'} \leqslant c_2 \| f \|_X \|g \|_{Y'}
\end {multline}
with some $c$, $c_1$ and $c_2$ independent of $f$ and $g$.

 %Observe that setting $X = Y = Z$ in Theorem~\ref {dpmxy} gives a slightly weakened version of Proposition~\ref {czobounda1p}.
Since $M$ is a positive operator and $M g \geqslant g$ almost everywhere for any locally summable $g$, conditions of Theorem~\ref {dpmxy}
imply that $X \subset Z$ and $Y' \subset Z'$ in the sense of continuous inclusions, which in turn implies that $X \subset Z \subset Y$.
Unlike the case $X = Y = Z$ it is presently unclear whether Theorem~\ref {dpmxy} admits a converse similar to
Theorem~\ref {a1necc} below.  In other words, if a suitable Calderon-Zygmund operator $T$ acts boundedly from $X$ to $Y$, does it follow that
there exists a lattice $Z$ satisfying the conditions of Theorem \ref {dpmxy}?

%We now pass to the proof of Theorem~\ref {dpmxy}.  By Theorem
%It is interesting to note that
There are still other approaches to establishing estimate \eqref {dpmxye}.
Let us now describe a recent result \cite {lerner2012}.  First, we need some preliminaries.
If $T$ is a Calderon-Zygmund operator with kernel $K$ then there is a maximal truncation operator
$$
T_\natural f (x) = \sup_{0 < \varepsilon < A} \left| \int_{\varepsilon < |y| < \nu} K (x, y) f (y)\, dy \right|, \quad x \in \mathbb R^n,
$$
associated with $T$ defined for all locally summable functions $f$.
It is well known (see, e.~g., \cite [Chapter~1, \S7] {stein1993}) that this operator is bounded in $\lclassg {p}$ for all $1 < p < \infty$.
The maximal truncated operator $T_\natural$ dominates the family of truncations
$$
T_{\varepsilon, A} f (x) = \int_{\varepsilon < |y| < \nu} K (x, y) f (y)\, dy, \quad x \in \mathbb R^n,
$$
of $T$, so this family of operators has a weak limit $T_*$ in $\lclassg {2}$ as $\varepsilon \to 0$ and $A \to \infty$,
and there exists some $a \in \lclassg {\infty}$ such that
$$
T f (x) = T_* f (x) + a (x) f (x)
$$
for all $f \in \lclassg {2}$ and almost all $x \in \mathbb R^n$.
Since multiplication by a bounded function $a$ is bounded in any lattice, boundedness of $T$ in a given lattice $X$ is thus implied by
boundedness of the maximal truncation operator~$T_\natural$ and vice versa.

A \emph {dyadic grid} $\mathcal D$
is a collection of cubes $Q$ in $\mathbb R^n$ with sides parallel to the coordinate axes such that
their lengths $\ell (Q)$ only take values $2^k$, $k \in \mathbb Z$, for any $Q, R \in \mathcal D$ we have
$Q \cap R \in \{Q, R, \emptyset\}$, and the cubes $\{Q \in \mathcal D \mid \ell (Q) = 2^k\}$ form a partition of
$\mathbb R^n$ for any $k \in \mathbb Z$.
A collection $\mathcal S = \{Q_j^k\} \subset \mathcal D$ is called a \emph {sparse family} of dyadic cubes
if it satisfies the following properties.
\begin {enumerate}
\item
Cubes $Q_j^k$ are pairwise disjoint in $j$ with $k$ fixed.
\item
If $\Omega_k = \bigcup_j Q_j^k$ then $\Omega_{k + 1} \subset \Omega_k$.
\item
$|\Omega_{k + 1} \cap Q_j^k| \leqslant \frac 1 2 |Q_j^k|$ for any $j$ and $k$.
\end {enumerate}
For any family of cubes $\mathcal S$ we define an operator
$$
\mathcal A_{\mathcal D, \mathcal S} f (x) = \mathcal A_{\mathcal S} f (x) = \sum_{Q \in \mathcal S} f_Q \chi_Q (x)
$$
acting on locally summable functions $f$, where as usual $f_Q = \frac 1 {|Q|} \int_Q f$.
It turns out that these operators with sparse families can be used to estimate Calderon-Zygmund operators
in the general setting.
\begin {theorem} [{\cite [Theorem~1.1] {lerner2012}}]
\label {dpm}
Suppose that $X$ is a Banach lattice on $\mathbb R^n$ having the Fatou property.
Then
$$
\|T_\natural f\|_X \leqslant c_{T, n} \sup_{\mathcal D, \mathcal S} \| \mathcal A_{\mathcal D, \mathcal S} |f| \|_X
$$
for any Calderon-Zygmund operator $T$ and locally summable function $f$ with compact support,
where the supremum is taken over arbitrary dyadic grids $\mathcal D$ and sparse families $\mathcal S \subset \mathcal D$.
\end {theorem}
Theorem~\ref {dpm} is based on a number of results that only recently were developed to a sufficient extent,
including a representation of Calderon-Zygmund operators as an average of dyadic shifts and
the local mean oscillation decomposition that represents every function $f \in S_0$ as $\mathcal A_{\mathcal S} f$ for some sparse family $\mathcal S$
in a given dyadic grid $\mathcal D$ with good pointwise control on $f - \mathcal A_{\mathcal S} f$;
see \cite {lerner2012} for a brief history of the techniques.

\begin {comment}
For the purposes of the present work
we can easily obtain the following conditions for boundedness.
\begin {theorem}
\label {dpmxy1}
Suppose that $X$, $Y$ and $Z$ are Banach lattices on $\mathbb R^n$ having the Fatou property,
the set $S_{00}$ of bounded functions with compact support is dense\footnote {Again,
if the norm of $X$ is order continuous then $S_{00}$ is dense in $X$
but not vice versa: take, for example, a variable exponent Lebesgue space $\lclassg {p (\cdot)}$ with $p = \infty$
on a certain bounded compact set $B$ and $1 \leqslant p < \infty$ elsewhere.}
in $X$, and the Hardy-Littlewood
maximal operator $M$ acts boundedly from $X$ to $Z$ and from $Y'$ to $Z'$.
Then any Calderon-Zygmund operator $T$ acts boundedly from $X$ to $Y$ and the same is true for the
maximal truncation operator $T_\natural$.
\end {theorem}
%Observe that setting $X = Y = Z$ in Theorem~\ref {dpmxy} gives a slightly weakened version of Proposition~\ref {czobounda1p}.
%Since $M$ is a positive operator and $M g \geqslant g$ almost everywhere for any locally summable $g$, conditions of Theorem~\ref {dpmxy}
%imply that $X \subset Z$ and $Y' \subset Z'$ in the sense of continuous inclusions, which implies that $X \subset Z \subset Y$.
%In contrast to the case $X = Y = Z$ it is presently unclear whether Theorem~\ref {dpmxy} admits a converse similar to
%Theorem~\ref {a1necc} below.  In other words, if a suitable Calderon-Zygmund operator $T$ acts boundedly from $X$ to $Y$, does it follow that
%there exists a lattice $Z$ satisfying the conditions of Theorem \ref {dpmxy}?
\end {comment}

As it was shown in \cite {lerner2012}, Theorem~\ref {dpm} has many interesting corollaries,
including the so-called $\apclass {2}$ conjecture and certain two-weight estimates.
Let us verify a replacement for the first line of the estimate \eqref {dpmxye} for
a suitable function $f \in X$ with operator $T_\natural$ in place of $T$.
By Theorem~\ref {dpm} there exists a dyadic grid $\mathcal D$ and a sparse family $\mathcal S \subset \mathcal D$ such that
$\|T_\natural f\|_Y \leqslant c \|\mathcal A_{\mathcal D, \mathcal S} |f| \|_Y$ with some constant $c$ independent of $f$.
Therefore there exist some $g \in Y'$, $\|g_0\|_{Y'} \leqslant 1$, such that
\begin {equation}
\label {dpm1}
\|T_\natural f\|_Y \leqslant 2 c \int \mathcal (A_{\mathcal D, \mathcal S} |f|)\, g.
\end {equation}
We may assume that $g \geqslant 0$.
The integral on the right-hand side of \eqref {dpm1} can now be %easily 
estimated using a kind of a stopping time argument \cite [(2.2)] {lerner2012}
which we are going to reproduce here.
Let $\{Q_j^k\} = \mathcal S$ and $\Omega_k$ be the cubes and sets in the definition of the sparse family $\mathcal S$,
and let $E_j^k = Q_j^k \setminus \Omega_{k + 1}$, so that $|E_j^k| \geqslant \frac 1 2 |Q_j^k|$ and $\{E_j^k\}$ is a collection of pairwise disjoint sets.
Then
\begin {multline*}
\int (\mathcal A_{\mathcal D, \mathcal S} |f|) \, g =
\sum_{j, k} \frac 1 {|Q_j^k|} \int_{Q_j^k} |f| \int_{Q_j^k} g =
\\
\sum_{j, k} |Q_j^k|\, \left(\frac 1 {|Q_j^k|} \int_{Q_j^k} |f|\right) \, \left(\frac 1 {|Q_j^k|} \int_{Q_j^k} g\right) \leqslant
\\
2 \sum_{j, k}  |E_j^k|\, \left(\frac 1 {|Q_j^k|} \int_{Q_j^k} |f|\right) \, \left(\frac 1 {|Q_j^k|} \int_{Q_j^k} g\right) =
\\
2 \sum_{j, k} \int \chi_{E_j^k}\, \left(\frac 1 {|Q_j^k|} \int_{Q_j^k} |f|\right) \, \left(\frac 1 {|Q_j^k|} \int_{Q_j^k} g\right) \leqslant
\\
2 \sum_{j, k} \int_{E_j^k} (M f) (M g) \leqslant
2 \int (M f) (M g),
\end {multline*}
which together with \eqref {dpm1} provides a suitable replacement for the first line in estimate \eqref {dpmxye}.
Another estimate for Calderon-Zygmund and certain other operators that can also be used to establish \eqref {dpmxye}
can be found in \cite {hytonen2012}.

 %Therefore we can continue estimate \eqref {dpm1} as
%$$
%\|T_\natural f\|_Y \leqslant 4 c \|M f\|_Z \|M g\|_{Z'} \leqslant C \|f\|_X \|g\|_{Y'} \leqslant C \|f\|_X
%$$
%using the assumptions of Theorem~\ref {dpmxy} with some constant $C$ independent of $f$.
%The proof of Theorem~\ref {dpmxy} is complete.

\section {Nondegenerate singular operators}

\label {nsos}

In this section we will try to give a more or less precise meaning to the nondegeneracy conditions that a singular operator
in Condition~2 of Theorem~\ref {themcr} must
satisfy in order to have a converse implication $2 \Rightarrow 1$ as well as
discuss certain restrictions on the spaces that are implied by boundedness of certain classes of operators.
Recall that the Muckenhoupt weights $\weightw \in \apclass {2}$ are precisely those for which the Hardy-Littlewood maximal operator is bounded in
the corresponding weighted space $\lclass {2} {\weightw^{-\frac 1 2}}$.
There is, however, a large class of operators that also characterize Muckenhoupt weights in this sense.
\begin {definition}
\label {nondego}
A mapping $T : \lclassg {2} \to \lclassg {2}$ is called $\apclass {2}$-nondegenerate with a constant $C$ if
boundedness of $T$ in a lattice $\lclass {2} {\weightw^{-\frac 1 2}}$ implies $\weightw \in \apclass {2}$ with constant $C$.
\end {definition}
This definition is stated for the general setting of a homogeneous space $S$ and measurable functions on $S \times \Omega$.
It is worth mentioning that for linear maps $T f (x, y) = [T_0 f (\cdot, y)] (x)$
that act in the first variable $x \in S$ only, i.~e. uniformly in $y \in \Omega$,
nondegeneracy of $T_0$ on $S$ implies nondegeneracy of $T$ on $S \times \Omega$; see \cite [Proposition~3.7] {rutsky2011en}.
For simplicity we will only work with a single variable from $S = \mathbb R^n$ in this section.

Although it is not clear yet how nondegeneracy in the sense of Definition~\ref {nondego} can be characterized in terms of the kernel of
a singular integral operator $T$, there are some useful sufficient conditions that illustrate this phenomenon.
 %Observe that a regular operator (such as convolution with a bounded summable function) cannot be $\apclass {2}$-nondegenerate since
\begin {definition}
\label {nondegso}
We say that a mapping $T : \lclassg {2} \to \lclassg {2}$ is nondegenerate along a direction $x_0 \in \mathbb R^n \setminus \{0\}$
if there exists a constant $c > 0$ such that for any ball $B \subset \mathbb R^n$ of radius $r > 0$ and any locally summable
nonnegative function $f$ supported on $B$
we have
\begin {equation}
\label {nondegsoeq}
|T f (x)| \geqslant c f_B
\end {equation}
for all $x \in B \pm r x_0$.
\end {definition}
It is well known that singularity of a mapping $T$ in the sense of Definition~\ref {nondegso} implies $\apclass {2}$-nondegeneracy of $T$;
in Proposition~\ref {nondegq} below we will establish a somewhat more general result.
In terms of the kernel $K$ of $T$ condition~\eqref {nondegsoeq} roughly means that $K (x, y)$ as a function of $|x - y|$
has to increase at $0$ as quickly and decay at infinity as slowly
as $|x - y|^{-n}$ along a certain direction; this statement is made more precise 
in Proposition~\ref {tndcs} below.
 %we will give a simple sufficient condition for nondegeneracy along a direction
 %of a Calderon-Zygmund operator.
It is not clear whether the class of mappings described by Definition~\ref {nondego}
is actually wider than that described by Definition~\ref {nondegso}.

Let $\mathcal S = \{Q_l\}$ be a collection of cubes or balls.
In addition to $\mathcal A_{\mathcal S}$ we introduce the following ``square'' averaging operator
$$
\mathcal A_{\mathcal S}^\square f (x) = \left(\sum_{Q \in \mathcal S} (f_Q)^2 \chi_Q (x)\right)^{\frac 1 2}
$$
for all locally summable functions $f$.  It is easy to see that if the cubes or balls from $\mathcal S$ are pairwise disjoint
then $\mathcal A_{\mathcal S}^\square f = \mathcal A_{\mathcal S} f$ almost everywhere for nonnegative functions $f$.
\begin {proposition}
\label {nondegq}
Suppose that a linear operator $T$ that is nondegenerate along a direction $x_0$
is bounded with norm $C$ in a Banach lattice $X$ having the Fatou property.
Then for any collection of cubes or balls $\mathcal S = \{Q_l\}$
we have
\begin {equation}
\label {asestimate}
\|\mathcal A_{\mathcal S}^\square f\|_X \leqslant c_a \left\|f \left(\sum_{l} \chi_{Q_l}\right)^{\frac 1 2} \right\|_X
\end {equation}
for all $f$ such that the right-hand part of \eqref {asestimate} is well-defined with
a constant $c_a$ independent of $f$ and $\mathcal S$.
\end {proposition}
To prove Proposition~\ref {nondegq} let $\mathcal S' = \{Q_l'\}$ with $Q_l' = Q_l + x_0$ being the cubes or balls $Q_l$
shifted by $x_0$
and set $f_l = f \chi_{Q_l}$.  We may assume that $f$ is nonnegative and that the right-hand part of \eqref {asestimate} is finite.
It follows that the sequence valued function $F = \{f_l\}$ belongs to $X (\lsclass {2})$ with
$\|F\|_{X (\lsclass {2})} = \left\|f \left(\sum_{l} \chi_{Q_l}\right)^{\frac 1 2} \right\|_X$.
Using the nondegeneracy assumption and the Grothendieck theorem (see, e.~g., \cite {krivine1973}) we can easily obtain the estimate
\begin {multline}
\label {nondegq1}
c^{-1} \left\|\left(\sum_l \chi_{Q_l'} (f_{Q_l})^2\right)^{\frac 1 2} \right\|_X \leqslant
\left\|\left(\sum_l \chi_{Q_l'} |T f_l|^2\right)^{\frac 1 2} \right\|_X \leqslant
\\
\|T F\|_{X (\lsclass {2})} \leqslant C K_G \|F\|_{X (\lsclass {2})} = C K_G \left\|f \left(\sum_{l} \chi_{Q_l}\right)^{\frac 1 2} \right\|_X,
\end {multline}
$K_G$ being the Grothendieck constant.  Repeating this estimate for function
$G = \{g_l\}$, $g_l = \chi_{Q_l'} f_{Q_l}$, in place of $F$ and with the order of $Q_l$ and $Q_l'$ reversed yields
\begin {multline}
\label {nondegq2}
c^{-1} \|\mathcal A_{\mathcal S}^\square f\|_X = c^{-1} \left\|\left(\sum_l \chi_{Q_l} (f_{Q_l})^2\right)^{\frac 1 2} \right\|_X \leqslant
\\
\left\|\left(\sum_l \chi_{Q_l} |T g_l|^2\right)^{\frac 1 2} \right\|_X \leqslant
\|T G\|_{X (\lsclass {2})} \leqslant C K_G \|G\|_{X (\lsclass {2})} =
\\
 %C K_G \left\|f \left(\sum_{l} \chi_{Q_l}\right)^{\frac 1 2} \right\|_X
%=
C K_G \left\|\left(\sum_l \chi_{Q_l'} (f_{Q_l})^2\right)^{\frac 1 2} \right\|_X.
\end {multline}
Combining \eqref {nondegq1} and \eqref {nondegq2} together yields \eqref {asestimate} with $c_a = (C c K_G)^2$.

The following corollary is essentially well known; see, e.~g., remarks after \cite [Lemma~5.2.2] {varpbook}.
\begin {corollary}
\label {nondegqc}
Suppose that a linear operator $T$ is nondegenerate along a direction $e$.
Then $T$ is $\apclass {2}$-nondegenerate.
\end {corollary}
Indeed, suppose that $T$ is bounded in $\lclass {2} {\weightw^{-\frac 1 2}}$ as in Definition~\ref {nondego}.
Taking in \eqref {asestimate}
a family $\mathcal S = \{B\}$ consisting of a single ball $B \subset \mathbb R^n$, $X = \lclass {2} {\weightw^{-\frac 1 2}}$
and a nonnegative locally summable function $f$ supported in $B$ yields
\begin {equation}
\label {nondegqc1}
f_B \left(\int_B \weightw \right)^{\frac 1 2} =
\|\mathcal A_{\mathcal S} f \|_X \leqslant c \left\|f\right\|_X = 
c \left(\int_B f^2 \weightw \right)^{\frac 1 2}.
\end {equation}
By rearranging the terms of \eqref {nondegqc1} we arrive at
$$
\left(f_B\right)^2 \leqslant c^2 \frac 1 {\int_B \weightw} \int_B f^2 \weightw,
$$
which is a well-known characterization of the Muckenhoupt weights from \cite [Chapter~5, \S1.4] {stein1993};
setting $f = (\weightw + \varepsilon)^{-1}$ and passing to the limit $\varepsilon \to 0$
quickly leads to \eqref {muckenhouptc}.

Observe that Proposition~\ref {nondegq} also implies that if a suitably nondegenerate operator $T$ acts boundedly in $X$ then
all operators $\mathcal A_{\mathcal S}$ with disjoint collections $\mathcal S$ of cubes or balls are uniformly bounded in $X$.
It is not clear in general how this property is related to other properties.  One, of course, immediately notices that
such operators $\mathcal A_{\mathcal S}$ are bounded in $\lclassg {p}$ for both $p = 1$ and $p = \infty$
so their uniform boundedness in a lattice $X$ does not imply per se that $X$ is $\apclass {1}$-regular.
However, and somewhat surprisingly, this implication holds true at least in the case of variable exponent Lebesgue spaces
if we also assume that $X$ is $p$-convex and $q$-concave for some $1 < p, q < \infty$; see \cite [Theorem~5.7.2] {varpbook}.
 %or \cite [Theorem~4.63] {varlsp}.
This rather involved result together with Proposition~\ref {nondegq}
provides at once the converse implication $3 \Rightarrow 1$ of Theorem~\ref {themcr} in the case of variable exponent Lebesgue spaces.
\begin {corollary}
\label {varexpa}
Suppose that $p (\cdot)$ is a measurable function on $\mathbb R^n$ such that
$1 < \essinf_{x \in \mathbb R^n} p (x) \leqslant \esssup_{x \in \mathbb R^n} p (x) < \infty$
and a linear operator $T$ is nondegenerate along a direction and bounded in
$\lclassg {p (\cdot)}$.  Then both $\lclassg {p (\cdot)}$ and $\lclassg {p' (\cdot)}$ are $\apclass {1}$-regular.
\end {corollary}
Thus not only is the converse to \cite [Theorem~5.39] {varlsp} true for nondegenerate operators,
which answers positively \cite [Problem~A.17] {varlsp},
there is also no need to involve the somewhat complicated machinery of the main results of this paper.

Now we give a standard condition sufficient for a Calderon-Zygmund operators to be nondegenerate along a direction.
\begin {proposition} [{\cite [Chapter~5, \S4.6] {stein1993}}]
\label {tndcs}
Suppose that $T$ is a Calderon-Zygmund operator with kernel $K$
and there exist some $u \in \mathbb R^n$ and a constant $c$ such that for any $x \in \mathbb R^n$ and $t \neq 0$ we have
\begin {equation}
\label {sdir}
|K (x, x + t u)| \geqslant c t^{-n}.
\end {equation}
Then $T$ is nondegenerate along the direction $x_0 = s u$ with some $s > 0$ and hence $T$ is $\apclass {2}$-nondegenerate.
\end {proposition}
The two typical examples are the Hilbert transform $H$ on $\mathbb R$ with kernel $K (x, y) = \frac {c_1} {x - y}$
and Riesz transforms $R_j$, $1 \leqslant j \leqslant n$ on $\mathbb R^n$
with kernels $K_j (x, y) = \frac {c_n (y_j - x_j)} {|y - x|^{n + 1}}$, where $c_n \neq 0$ are some constants.
It is evident that these kernels satisfy condition~\eqref {sdir} for $u = e_j$, $e_j$ being the $j$-th coordinate basis vector of $\mathbb R^n$.

For completeness, let us prove Proposition~\ref {tndcs}.
Indeed, condition \eqref {kernelest} on the kernel $K$ implies that
\begin {equation}
\label {k1ndc}
|K (x, x + r [-x_0 + v]) - K (x, x - r x_0)| \leqslant C_K \frac {|r v|^\gamma} {|r x_0|^{n + \gamma}}
\leqslant c'  s^{-n - \gamma} r^{-n}
\end {equation}
for all $v \in \mathbb R^n$, $|v| < \frac 1 2 |x_0| = \frac 1 2 s |u|$,
and any $r \neq 0$ with some constant $c'$ independent of $s$.
By taking $s$ sufficiently large we may assume that \eqref {k1ndc} holds true for all $|v| < 1$.
Therefore \eqref {sdir} implies that
\begin {multline*}
|T f (x)| =
\\
\left|
K (x, x - r x_0) \int\limits_B f (y) dy + \int\limits_B \left[ K (x, y) - K (x, x - r x_0)\right] f (y) dy
\right| \\
\geqslant |K (x, x - r x_0)| \int\limits_B f (y) dy -
\left|\int\limits_B \left[ K (x, y) - K (x, x - r x_0)\right] f (y) dy\right|
\\
\geqslant |K (x, x - r x_0)| \int\limits_B f (y) dy -
\int\limits_B \left| K (x, y) - K (x, x - r x_0)\right| |f (y)| dy
\\
\geqslant \int\limits_B f (y) dy \cdot \left( c (r s)^{-n} - c' s^{-n - \gamma} r^{-n}\right) =
(c s^{-n} - c' s^{-n - \gamma}) f_B
\end {multline*}
for any $f$ supported on the ball $B \subset \mathbb R^n$ having radius $r$ and centered at the origin and for any $x \in B \pm r x_0$.
Choosing $s$ sufficiently large yields \eqref {nondegsoeq}, so $T$ is indeed nondegenerate along the direction $x_0$
and is therefore $\apclass {2}$-nondegenerate by Corollary~\ref {nondegqc}.  The proof of Proposition~\ref {tndcs} is complete.

\end {unabridged}

\section {A lemma about $\apclass {p}$-regularity}

\label {aprmainlemma}

\begin {unabridged}
In this section we establish the following auxiliary result that we will need in Section~\ref {noaro} below.
\end {unabridged}
\begin {theorem}
\label {a1apt}
Suppose that $X$ is a Banach lattice of measurable functions on $(S \times \Omega, \mu \times \nu)$
such that $X$ satisfies the Fatou property, and
\begin {enumerate}
\item
$X$ is $\apclass {p}$-regular with constants $(c_1, m_1)$
for some $1 < p < \infty$, %and
\item
$X^\delta$ is $\apclass {1}$-regular with constants $(c_2, m_2)$ for
some $\delta > 0$.
\end {enumerate}
Then lattice $X$ is $\apclass {1}$-regular with an estimate for the
constants depending only on the corresponding $\apclass {p}$-regularity constants of $X$,
$\apclass {1}$-regularity constants of $X^\delta$ and the value of $\delta$.
\end {theorem}
This theorem is easily derived from the corresponding result for $\apclass {p}$ weights with the help of a fixed point argument.
\begin {lemma}
\label {a1apl}
Suppose that a weight $\weightw$ on $(S \times \Omega, \mu \times \nu)$ satisfies
$\weightw \in \apclass {p}$ and $\weightw^\delta \in \apclass {1}$ with some $1 < p < \infty$ and $\delta > 0$.
Then $\weightw \in \apclass {1}$ with an estimate for the constants depending only on $\delta$, the corresponding constants
of the $\apclass {p}$ condition for $\weightw$ and the $\apclass {1}$ condition for $\weightw^\delta$.
\end {lemma}
Lemma~\ref {a1apl} is essentially a particular case $X = \lclass {\infty} {\weightw}$
of Theorem~\ref {a1apt}.
\begin {unabridged}
This result is suggested by a very simple observation: by the factorization of $\apclass {p}$ weights
(see, e.~g., \cite [Chapter~5, \S5.3] {stein1993}) we have $\weightw = \omega_0 \omega_1^{1 - p}$ with $\omega_j \in \apclass {1}$,
and since we also have $\weightw^\delta \in \apclass {1}$, $\weightw$ is bounded away from $0$ on every ball, which indicates that
the singularities of the denominator factor $\omega_1$ have to be dominated by the singularities of the nominator factor
$\omega_0$ in some sense
and $\omega_1$ should essentially
cancel out in this factorization.
\end {unabridged}
To prove Lemma~\ref {a1apl},
fix some $\omega \in \Omega$ such that $\weightw (\cdot, \omega) \in \apclass {p}$ and $\weightw^\delta (\cdot, \omega) \in \apclass {1}$,
and let $B (x, r) \subset S$, $x \in S$, $r > 0$, be an arbitrary ball of $S$.
Then sequential application of the $\apclass {p}$ condition satisfied by weight $\weightw$,
the Jensen inequality with convex function $t \mapsto t^{-\delta (p - 1)}$, $t > 0$, and the $\apclass {1}$ condition satisfied by the weight $\weightw^\delta$
yields
\begin {multline}
\label {a1condb}
\frac 1 {\nu (B (x, r))} \int_{B (x, r)} \weightw (u, \omega) d\nu (u) \leqslant
\\
c \left[ \frac 1 {\nu (B (x, r))} \int_{B (x, r)} [\weightw (u, \omega)]^{-\frac 1 {p - 1}} d\nu (u) \right]^{-(p - 1)} =
\\
c \left[ \frac 1 {\nu (B (x, r))} \int_{B (x, r)} [\weightw (u, \omega)]^{-\frac 1 {p - 1}} d\nu (u) \right]^{-\delta (p - 1) \cdot \frac 1 \delta} \leqslant
\\
c \left[ \frac 1 {\nu (B (x, r))} \int_{B (x, r)} [\weightw (u, \omega)]^\delta d\nu (u) \right]^{\frac 1 \delta} \leqslant
c' \weightw (x, \omega)
\end {multline}
for almost all $x \in S$ with some constants $c$ and $c'$ depending only on the corresponding constants
of the $\apclass {p}$ condition for $\weightw$, the $\apclass {1}$ condition for $\weightw^\delta$ and the value of $\delta$.
Since $\omega$, $x$ and $B$ are arbitrary, \eqref {a1condb} implies that $\weightw \in \apclass {1}$ with the necessary
estimates of the constants, which concludes the proof of Lemma~\ref {a1apl}.

In order to reduce Theorem~\ref {a1apt} to Lemma~\ref {a1apl} we need to show that under the conditions of Theorem~\ref {a1apt}
an arbitrary function
$f \in X$ has a majorant $\weightw$ such that with the appropriate estimates on the constants
$\weightw$ is an $\apclass {p}$-majorant of $f$ in $X$
and simultaneously $\weightw^\delta$ is an $\apclass {1}$-majorant of $|f|^\delta$
in $X^\delta$.  At a first glance it may seem that there is little reason to suspect existence of a common majorant in
sets that look vastly different (for example,
a majorant $\weightw$ such that $\weightw^\delta \in \apclass {1}$ may not even be locally summable in the first variable,
while on the other hand a majorant $\weightw \in \apclass {p}$ may vanish near some points);
however, careful application of the celebrated Ky-Fan--Kakutani
fixed point theorem allows us to establish the existence of a common majorant in this setting with relative ease.
\begin {theoremnn} [{\cite {fanky1952}}]
Suppose that $K$ is a compact set in a locally convex linear topological space.
Let $\Phi$ be a mapping from $K$ to the set of nonempty convex compact subsets of $K$.
If the graph
$$
\Gamma (\Phi) = \{ (x, y) \in K \times K \mid y \in \Phi (x) \}
$$
of $\Phi$ is closed in $K \times K$
then $\Phi$ has a fixed point, i.~e. $x \in \Phi (x)$ for some $x \in K$.
\end {theoremnn}
We will also need the following sets of nonnegative a. e. measurable functions $\weightw$ on
$(S \times \Omega, \nu \times \mu)$ (see also \cite [Section~3] {rutsky2011en}):
$$
\abp {p} {C} = \left\{ \weightw \mid \esssup_{\omega \in \Omega} \|M\|_{\lclass {p} {\weightw^{-\frac 1 p} (\cdot, \omega)}} \leqslant C \right\};
$$
$$
\abp {1} {C} = \left\{ \weightw \mid \esssup \frac {M \weightw} {\weightw} \leqslant C \right\}.
$$
These are the sets of Muckenhoupt weights with fixed bounds on the constants
(``the Ball of $\apclass {p}$'').
 %Since by our conventions $0 \in \abp {p} {C}$ these sets are nonempty for all $C \geqslant 0$.
\begin {proposition} [{\cite [Proposition~3.4] {rutsky2011en}; see also \cite[Lemma 4.2]{kalton1994}}]
\label {bmoconv}
Suppose that $1 \leqslant p < \infty$ a. e. and $C \geqslant 0$.
The set $\abp {p} {C}$ is a nonempty convex cone which is also logarithmically convex and closed in measure.
\end {proposition}
\begin {unabridged}
The proof of convexity and closedness is quite routine; the logarithmic convexity is a bit harder but we will not need it under
the assumptions of Theorem~\ref {a1apt}.
\end {unabridged}

We are now ready to prove Theorem~\ref {a1apt}.
The technical details of this proof as well as the general pattern are similar to the main result of
\cite {rutsky2011en}.
By using \cite [Proposition~3.6] {rutsky2011en}
it is sufficient to establish the existence of a suitable majorant
for every function $f \in X$, $\|f\|_X \leqslant 1$, such that $E = \supp f$ has positive finite measure and $f \geqslant \beta$ on $E$
with some $\beta > 0$,
since the set of such functions is dense in measure in the nonnegative part of the closed unit ball $B$ of $X$.
We fix such a function~$f$.

 %Now we only have the inclusion $X \subset \lclass {1} {a}$ for some $a \in X'$, $a > 0$ a.~e.,
%so the set $B$ defined as before may not necessarily be compact.
%However, simple estimates show that
By Proposition~\ref {xpos} there exists some function $a \in X'$, $\|a\|_{X'} = 1$, such that $a > 0$ almost everywhere.
This implies that for any $u \in B$ we have
$\int |u| a \leqslant \|u\|_X \|a\|_{X'} \leqslant 1$, i.~e. $\|u\|_{\lclass {1} {a^{-1}}} \leqslant 1$.
Let $0 < \alpha \leqslant \beta \leqslant 1$ be a sufficiently small number to be determined later, and let
$$
D = \{  \chi_E \log g \mid g \in B, \, g \geqslant \chi_E \alpha \}.
$$
It is easy to see that $D$ is a bounded set in $Y = \lclass {2} {a^{-\frac 1 2}}$ for any given $E$ and $\alpha$ because
$$
\int_{E \cap \{g < 1\}} |\log g|^2 a \leqslant |\log \alpha|^2 \|\chi_E\|_X \|a\|_{X'} \leqslant \frac 1 \beta |\log \alpha|^2
$$
and
$$\int_{E \cap \{g \geqslant 1\}} |\log g|^2 a =
\int_{E \cap \{g \geqslant 1\}} 4 \left|\log \left(g^{\frac 1 2}\right)\right|^2 a \leqslant 4 \int |g| a \leqslant 4
$$
for any $\chi_E \log g \in D$; $D$ is convex because $B$ is logarithmically convex and $D$ is closed in measure,
so $D$ is compact in the weak topology of~$Y$.

Observe that since $\apclass {1}$-regularity of $X$ implies $\apclass {1}$-regularity of $X^\gamma$ for all $0 < \gamma < 1$ we may assume that
$0 < \delta < 1$, otherwise the conclusion of Theorem~\ref {a1apt} is immediate.
 %It is easy to verify
%in the same way as in Proposition~\ref {bmoconv} that the set
%$\left[\abp {1} {C}\right]^{\frac 1 \delta} = \{\weightw \mid \weightw^\delta \in \abp {1} {C} \}$ is convex using the fact that
%$(a + b)^\delta \leqslant a^\delta + b^\delta$ for all numbers $a, b \geqslant 0$, and closedness in measure of
%$\left[\abp {1} {C}\right]^{\frac 1 \delta}$ follows immediately from Proposition~\ref {bmoconv}.
We define a set-valued map $\Phi$ in~$D \times D$ %acting into subsets of $D \times D$
~by
\begin {multline*}
\Phi ((\log \weightu, \log \weightv)) = \left\{ (\log \weightu_1, \log \weightv_1) \mid
\weightu_1, \weightv_1 \in X,
\right.
\\
\left.
\weightu_1 \in B \cap \abp {p} {c_1},\,\, \weightv_1^\delta \in B \cap \abp {1} {c_2},
\right.
\\
\left.
f \vee (\weightu \vee \weightv) \leqslant A \left(\weightu_1 \wedge \weightv_1\right) \right\}
\end {multline*}
Since for any $(\log \weightu, \log \weightv) \in D \times D$ we have
$\weightw = f \vee \weightu \vee \weightv \in X$ with
$\|\weightw\|_X \leqslant 3$ and by the assumptions there exist some $a, b \in X$ such that
$a \in \abp {p} {c_1}$, $b^\delta \in \abp {1} {c_2}$, $a \geqslant \weightw$, $b \geqslant \weightw$ and
$\|a\|_X \leqslant 3 m_1$, $\|b\|_X \leqslant (3 m_2)^\delta$.
Thus choosing $A = (3 m_1) \vee (3 m_2)^{\frac 1 \delta}$ and $\alpha = \beta \wedge A^{-1}$
yields $(\log \weightu_1, \log \weightv_1) \in \Phi ((\log \weightu, \log \weightv))$ with
$\weightu_1 = \frac 1 A a$ and $\weightv_1 = \frac 1 A b$, so $\Phi$ takes nonempty values.
 %Nonemptiness of the values of $\Phi$ is implied by the $\apclass {p}$-regularity of $X$
%and $\apclass {1}$-regularity of $X^\delta$ by the assumptions of Theorem~\ref {a1apt}.
The condition $f \vee (\weightu \vee \weightv) \leqslant A \left(\weightu_1 \wedge \weightv_1\right)$ is
of course equivalent to (and a shorthand for) the six inequalities
$f \leqslant A \weightu_1$, $f \leqslant A \weightv_1$,
$\weightu \leqslant A \weightu_1$, $\weightv \leqslant A \weightu_1$,
$\weightu \leqslant A \weightv_1$ and $\weightv \leqslant A \weightv_1$.
It is easy to see using Proposition~\ref {bmoconv} %and remarks after it
that the graph $\Gamma$ of $\Phi$ is a convex set and $\Gamma$ is closed with respect to the convergence in measure.
Let us verify that $\Gamma$ is closed in $Y \times Y$.
Indeed, the weak topology of $Y \times Y$ is metrizable on a bounded set $D \times D$.
If $x_j \in \Gamma$ and $x_j \to x \in Y \times Y$ then there exists some sequence $y_j$ of convex combinations of $x_j$ such that
$y_j \to x$ in the strong topology of $Y \times Y$, and $y_j \in \Gamma$ by the convexity of $\Gamma$.
Strong convergence in $Y$ implies convergence
in measure, so $y_j \to x$ in measure.  Since $\Gamma$ is closed in measure, it follows that $x \in \Gamma$ and thus $\Gamma$ is indeed closed in
$Y \times Y$.  From this we also infer that the values of $\Phi$ are convex and closed in the compact set $D \times D$ and thus
they are compact in $Y \times Y$.

By the Ky Fan--Kakutani fixed point theorem there exists some $(\log \weightu, \log \weightv) \in D \times D$ such that
$(\log \weightu, \log \weightv) \in \Phi ((\log \weightu, \log \weightv))$.
This implies that $\weightu$ and $\weightv$ are pointwise equivalent to one another with the constant of equivalence depending only on $A$
(which, in turn, only depends on the values of $m_1$, $m_2$ and $\delta$),
and so $\weightw = A \weightu$ is a majorant of $f$ such that $\weightw \in \apclass {p}$ and $\weightw^\delta \in \apclass {1}$ with the
appropriate estimates on the constants.  By Lemma~\ref {a1apl} it follows that $\weightw \in \apclass {1}$ with suitable estimates on the constants,
which concludes the proof of Theorem~\ref {a1apt}.

We will need the following proposition, which is a simple consequence of duality and the properties of $\apclass {p}$ weights.
\begin {proposition} [{\cite [Proposition~2.3] {rutsky2011en}}]
\label {aptoconj}
Suppose that $X$ is a Banach lattice on $(S \times \Omega, \nu \times \mu)$
such that $X'$ is a norming space for $X$.
If $X'$ is $\apclass {1}$-regular then $X^{\frac 1 q}$ is $\apclass {1}$-regular for all $q > 1$.
If $X'$ is $\apclass {p}$-regular with some $p > 1$ then
$X^{\frac 1 p}$ is $\apclass {1}$-regular.
\end {proposition}

Theorem~\ref {a1apt} has an interesting immediate application.
\begin {proposition}
\label {ainfainf}
Let $X$ be a Banach lattice on $(S \times \Omega, \nu \times \mu)$ having the Fatou property.
Suppose that both $X$ and $X'$ are $\apclass {\infty}$-regular.
Then both $X$ and $X'$ are $\apclass {1}$-regular.
\end {proposition}
Indeed, since $X$ and $X'$ are $\apclass {\infty}$-regular, they are also $\apclass {p}$-regular with some $p > 1$, which by Proposition~\ref {aptoconj}
means that both $X'^{\frac 1 p}$ and $X^{\frac 1 p}$ are $\apclass {1}$-regular, and it remains to apply Theorem~\ref {a1apt} to $X$ and $X'$
with $\delta = \frac 1 p$.

\begin {corollary}
\label {czobounda1pc}
Suppose that $X$ is a Banach lattice on $\mathbb R^n$ having the Fatou property,
and both $X$ and $X'$ are $\apclass {\infty}$-regular.
Then any Calderon-Zygmund operator $T$ is bounded in~$X$.
\end {corollary}
This corollary, which strengthens Proposition~\ref {czobounda1p},
immediately follows from Proposition~\ref {ainfainf} and Propositon~\ref {czobounda1p}.

\section {Necessity of $\apclass {1}$-regularity}

\label {noaro}

\begin {unabridged}
In this section we establish the converse implication $3 \Rightarrow 1$ of Theorem~\ref {themcr}.
We will need the following fairly well known result, the proof of which %along with a brief history
in the present setting can be found in \cite [Theorem~2.6] {rutsky2011en}.
\end {unabridged}
\begin {abridged}
The proof of the following result can be found in \cite [Theorem~2.6] {rutsky2011en}.
\end {abridged}
\begin {theorem}
\label {btsb}
Suppose that $Y$ is a Banach lattice on $(S \times \Omega, \nu \times \mu)$ with an order continuous norm.
If a linear operator $T$ is bounded in $Y^{\frac 1 2}$
then for every $f \in Y'$, $m > 1$ and $a > K_G \|T\|_{Y^{\frac 1 2}}$, $K_G$ being the Grothendieck constant,
there exists a majorant
$\weightw \geqslant |f|$, $\|\weightw\|_{Y'} \leqslant \frac m {m - 1} \|f\|_{Y'}$, such that
$\|T\|_{\lclass {2} {\weightw^{-\frac 1 2}} \to \lclass {2} {\weightw^{-\frac 1 2}}} \leqslant a \sqrt m$.
\end {theorem}
\begin {unabridged}
This theorem essentially says that for suitably nondegenerate operators $T$ boundedness of $T$
in a lattice $Y^{\frac 1 2}$ implies that $Y'$ is $\apclass {2}$-regular,
which binds the boundedness property of certain operators in a lattice back to a regularity property for some related lattices.
The proof of Theorem~\ref {btsb} given in \cite [\S6] {rutsky2011en} is merely a slight refinement of the proof of
\cite[Theorem 3.5]{kisliakov1999},
which is in turn a variant of the well-known Maurey--Krivine factorization theorem (see \cite {cbs}).
For the first time these ideas were exploited in a similar context in \cite {rubiodefrancia1987}.
\end {unabridged}

\begin {theorem} [{\cite [Theorem~1.6] {rutsky2011en}}]
\label {bmomdiv}
Suppose that $X$ is a Banach lattice on $(S \times \Omega, \nu \times \mu)$
having the Fatou property.  Suppose also that
$X \lclassg {q}$ for some $1 < q < \infty$ is a Banach lattice and $X \lclassg {q}$ is $\apclass {p}$-regular for some $1 \leqslant p < \infty$.
Then $X$ is $\apclass {p + 1}$-regular.
\end {theorem}
\begin {unabridged}
Theorem~\ref {bmomdiv}, which is rather involved, is a direct precursor
to a very deep and nontrivial fact that the so-called BMO-regularity property is self-dual at least
for Banach lattices having the Fatou property; see \cite {kalton1994}, \cite {kisliakov2002en}, \cite {rutsky2011en}.
\end {unabridged}

\begin {theorem}
\label {a1necc}
Suppose that $X$ is a Banach lattice of measurable functions on $(S \times \Omega, \nu \times \mu)$ 
such that $X$ is $p$-convex and $q$-concave for some $1 < p, q < \infty$ and $X$ satisfies the Fatou property.
Let $T$ be a linear operator on $\lclass {2} {S \times \Omega}$ such that both $T$ and $T^*$ are $\apclass {2}$-nondegenerate
and $T$ acts boundedly in $X$ and in all $\lclassg {s}$ for $1 < s < \infty$.  Then lattices $X$ and $X'$ are $\apclass {1}$-regular.
\end {theorem}
Let us now prove Theorem~\ref {a1necc}.
By the $p$-convexity condition $X^p$ is also a Banach lattice with the Fatou property, and so
 %$X_{\phantom t}^{p (1 - \theta)} \lclassg {t}^\theta$
$X^{p (1 - \theta)} \lclassg {t}^\theta$
is also a Banach lattice for all $1 \leqslant t \leqslant \infty$ and $0 < \theta < 1$.
Choosing $\theta = 1 - \frac 1 p$ shows that $Y_s = X \lclassg {s}$ is a Banach lattice for all sufficiently large $s$.
Lattice $Y_s$ satisfies the Fatou property
and has order continuous norm (because $\lclassg {s}$ has order continuous norm for $s < \infty$).
Since $T$ is bounded in $X$ and in $\lclassg {s}$ for all $1 < s < \infty$, by the interpolation theorem
mentioned in Section~\ref {preliminaries} operator $T$ is also bounded in $X^{\frac 1 2}_{\phantom s} \lclassg {s}^{\frac 1 2} = Y_{s}^{\frac 1 2}$
for all $1 < s < \infty$.  Theorem~\ref {btsb} and $\apclass {2}$-nondegeneracy of $T$ then imply that
lattice $Y_{s}' = X' \lclassg {s'}$ is $\apclass {2}$-regular for all sufficiently large~$s$.
By Theorem~\ref {bmomdiv} it follows that lattice
$X'$ is $\apclass {3}$-regular, and furthermore by Proposition~\ref {aptoconj} lattice $X^{\frac 1 3}$ is $\apclass {1}$-regular.
Since the convexity assumptions of Theorem~\ref {a1necc} imply that lattices $X$ and $X'$ have order continuous norm,
we have $X' = X^*$ and $X = (X')^*$, and moreover $X \cap \lclassg {2}$ is dense in $X$ and $X' \cap \lclassg {2}$ is dense in $X'$,
so the duality relation $\int (T f) g = \int f (T^* g)$ for $f \in X \cap \lclassg {2}$ and $g \in X' \cap \lclassg {2}$
shows that
boundedness of $T$ in $X$ implies boundedness of the conjugate operator $T^*$ in $X'$ and vice versa.
Repeating the argument above with lattice $X'$ in place of $X$ (which is a $q'$-convex lattice since $X$ is $q$-concave) and
operator $T^*$ in place of $T$ shows that
lattice $X$ is $\apclass {3}$-regular and lattice $(X')^{\frac 1 3}$ is $\apclass {1}$-regular.
Finally, we apply Theorem~\ref {a1apt} to $X$ and to $X'$ with $p = 3$ and $\delta = \frac 1 3$, which establishes that
lattices $X$ and $X'$ are both $\apclass {1}$-regular. The proof of Theorem~\ref {a1necc} is complete.

%Proposition~\ref {czobounda1p} together with Theorem~\ref {a1necc} yield at once the following property of $\apclass {p}$-regular lattices.
%\begin {corollary}
%Suppose that $X$ is a Banach lattice of measurable functions on $(\mathbb R^n \times \Omega, \nu \times \mu)$ 
%such that $X$ is $p$-convex and $q$-concave for some $1 < p, q < \infty$ and $X$ satisfies the Fatou property.
%Suppose also that lattice $X$ is $
%\end {corollary}

\begin {abridged}

\section {Concluding remarks}

\label {concrems}

The theory of Calderon-Zygmund operators naturally generalizes to the spaces of homogeneous type.
However, we cannot just replace $\mathbb R^n$ by a space of homogeneous type in the statement of Theorem~\ref {themcr}
because it is not clear for which of these spaces %$\Omega$
there is at least one suitably nondegenerate linear operator $R$.
 %It is easy to see that the proof of Theorem~\ref {themcr} also works in the vector-valued case,
%i.~e. for lattices of measurable functions like $X (\lsclass {r})$, where $X$ is a lattice on $\mathbb R^n$.

The $p$-convexity and $q$-concavity assumptions of Theorem~\ref {themcr} are probably superfluous
(but so far we can only say for sure that they are not used in the implication $1 \Rightarrow 2$).
I conjecture that these assumptions are actually a consequence of any of the conditions
of Theorem~\ref {themcr}; that Condition~1 implies $p$-convexity and $q$-concavity with some $1 < p, q < \infty$
is known to hold true at least in the case of the variable exponent Lebesgue spaces (see, e.~g., \cite [Theorem~4.7.1] {varpbook}),
and it seems that it is possible to adapt the same argument to cover suitable nondegenerate singular integral operators as well.
Recently in \cite [Theorem~5.42] {varlsp} it was established that if all Riesz transforms $R_j$ are bounded in $\lclassg {p (\cdot)}$ then
the exponent $p (\cdot)$ is bounded away from $1$ and $\infty$.

Furthermore, it is easy to see that the assumptions of $p$-convexity and $q$-concavity could be eliminated from Theorem~\ref {a1necc} if
we assume boundedness of $T$ in $X^{\frac 1 2}$ and $(X')^{\frac 1 2}$ instead of just $X$.
 %and get rid of $p$-convexity assumption in Theorem~\ref {a1apt} as described at the end of this section.
It is, however, unclear whether boundedness of a Calderon-Zygmund operator $T$ in $X$ implies
its boundedness in $X^{\frac 1 2}$; it seems plausible because $T$ acts boundedly from $\lclassg {\infty}$ to $\BMO$,
but as far as I know,
all available interpolation results that make it possible to replace $\lclassg {\infty}$ by $\BMO$ as an endpoint in the appropriate interpolation scale
(see, e.~g., \cite {kopaliani2009}, \cite {rutsky2013iv1})
work only under the assumption that both $X^\alpha$ and $(X')^\alpha$ are $\apclass {1}$-regular
for some $\alpha > 0$, which is part of what we are trying to establish in this setting.
 %which are, however, true for all linear operators and not just Calderon-Zygmund,

There is a different approach to the implication $3 \Rightarrow 1$ of Theorem~\ref {themcr} that works at least in certain cases.
Let $\mathcal S = \{Q_l\}$ be a collection of cubes or balls.
We define operators
$$
\mathcal A_{\mathcal S} f (x) = \sum_{Q \in \mathcal S} f_Q \chi_Q (x)
$$
and
$$
\mathcal A_{\mathcal S}^\square f (x) = \left(\sum_{Q \in \mathcal S} (f_Q)^2 \chi_Q (x)\right)^{\frac 1 2}
$$
for all locally summable functions $f$; here $f_Q = \frac 1 {|Q|} \int_Q f (z) dz$ for $Q \in \mathcal S$.
It is easy to see that if the cubes or balls from $\mathcal S$ are pairwise disjoint
then $\mathcal A_{\mathcal S}^\square f = \mathcal A_{\mathcal S} f$ almost everywhere for nonnegative functions $f$.
\begin {proposition}
\label {nondegq}
Suppose that a singular integral operator $T$ with kernel $K$ satisfying \eqref {sdir}
is bounded with norm $C$ in a Banach lattice $X$ having the Fatou property.
Then for any collection of cubes or balls $\mathcal S = \{Q_l\}$
we have
\begin {equation}
\label {asestimate}
\|\mathcal A_{\mathcal S}^\square f\|_X \leqslant c_a \left\|f \left(\sum_{l} \chi_{Q_l}\right)^{\frac 1 2} \right\|_X
\end {equation}
for all $f$ such that the right-hand part of \eqref {asestimate} is well-defined with
a constant $c_a$ independent of $f$ and $\mathcal S$.
\end {proposition}
Observe that Proposition~\ref {nondegq} also implies that if a suitably nondegenerate operator $T$ acts boundedly in $X$ then
all operators $\mathcal A_{\mathcal S}$ with disjoint collections $\mathcal S$ of cubes or balls are uniformly bounded in $X$.
In particular, it is well known that taking collections $\mathcal S$ consisting of single cubes implies that if $X = \lclass {p} {\weightw^{-\frac 1 p}}$,
$1 < p < \infty$,
then $\weightw \in \apclass {p}$ and thus $X$ is $\apclass {1}$-regular.
It is not clear in general whether either \eqref {asestimate} or uniform boundedness of $\mathcal A_{\mathcal S}$ is related to other properties of interest.
One, of course, immediately observes that
such operators $\mathcal A_{\mathcal S}$ are bounded in $\lclassg {p}$ for both $p = 1$ and $p = \infty$
so their uniform boundedness in a lattice $X$ does not necessarily mean that $X$ is $\apclass {1}$-regular.
However, and somewhat surprisingly, this implication holds true at least in the case of variable exponent Lebesgue spaces $X = \lclassg {p (\cdot)}$
if we also assume that $X$ is $p$-convex and $q$-concave for some $1 < p, q < \infty$; see, e.~g.,
\cite [Theorem~5.7.2] {varpbook}. % or \cite [Theorem~4.63] {varlsp}.
Thus not only the converse to \cite [Theorem~5.39] {varlsp} is true for nondegenerate operators,
which answers positively \cite [Problem~A.17] {varlsp},
but there is also no need to involve the complicated machinery of the main results of the present work.

Let us prove Proposition~\ref {nondegq}. First, observe that \eqref {sdir} implies by \cite [Chapter~5, \S4.6] {stein1993} that
there exists a constant $c > 0$ and some $x_0 \in \mathbb R^n \setminus \{0\}$
such that for any ball $B \subset \mathbb R^n$ of radius $r > 0$ and any locally summable
nonnegative function $f$ supported on $B$
we have
\begin {equation}
\label {nondegsoeq}
|T f (x)| \geqslant c f_B
\end {equation}
for all $x \in B \pm r x_0$.
Let $\mathcal S' = \{Q_l'\}$ with $Q_l' = Q_l + x_0$ being the cubes or balls $Q_l$
shifted by $x_0$
and set $f_l = f \chi_{Q_l}$.  We may assume that $f$ is nonnegative and that the right-hand part of \eqref {asestimate} is finite.
It follows that the sequence valued function $F = \{f_l\}$ belongs to $X (\lsclass {2})$ with
$\|F\|_{X (\lsclass {2})} = \left\|f \left(\sum_{l} \chi_{Q_l}\right)^{\frac 1 2} \right\|_X$.
Using the nondegeneracy assumption in the form of \eqref {nondegsoeq}
and the Grothendieck theorem (see, e.~g., \cite {krivine1973}) we can easily obtain the estimate
\begin {multline}
\label {nondegq1}
c^{-1} \left\|\left(\sum_l \chi_{Q_l'} (f_{Q_l})^2\right)^{\frac 1 2} \right\|_X \leqslant
\left\|\left(\sum_l \chi_{Q_l'} |T f_l|^2\right)^{\frac 1 2} \right\|_X \leqslant
\\
\|T F\|_{X (\lsclass {2})} \leqslant C K_G \|F\|_{X (\lsclass {2})} = C K_G \left\|f \left(\sum_{l} \chi_{Q_l}\right)^{\frac 1 2} \right\|_X,
\end {multline}
$K_G$ being the Grothendieck constant.
On the other hand, repeating this estimate for function
$G = \{g_l\}$, $g_l = \chi_{Q_l'} f_{Q_l}$, in place of $F$ and with the order of $Q_l$ and $Q_l'$ reversed
shows that
\begin {multline}
\label {nondegq2}
c^{-1} \|\mathcal A_{\mathcal S}^\square f\|_X = c^{-1} \left\|\left(\sum_l \chi_{Q_l} (f_{Q_l})^2\right)^{\frac 1 2} \right\|_X \leqslant
\\
\left\|\left(\sum_l \chi_{Q_l} |T g_l|^2\right)^{\frac 1 2} \right\|_X \leqslant
\|T G\|_{X (\lsclass {2})} \leqslant C K_G \|G\|_{X (\lsclass {2})} =
\\
C K_G \left\|\left(\sum_l \chi_{Q_l'} (f_{Q_l})^2\right)^{\frac 1 2} \right\|_X.
\end {multline}
Combining \eqref {nondegq1} and \eqref {nondegq2} together yields \eqref {asestimate} with $c_a = (C c K_G)^2$.

There is an interesting generalization of Proposition~\ref {czobounda1p}
which is easily obtained from certain less classical results. % that we now briefly describe.
It is well known (see, e.~g., \cite {alvarezperez1994}, \cite {jawerthtorchinsky1985}) that
\begin {equation}
\label {flsharpczo}
M^\sharp_\lambda (T f) \leqslant c M f
\end {equation}
almost everywhere for a wide variety of operators $T$ including Calderon-Zygmund operators
and
all locally summable functions $f$ with $c$ independent of $f$,
where $M^\sharp_\lambda$ is the Str\"omberg local sharp maximal function.
$S_0$ denotes the set of all measurable functions $f$ on $\mathbb R^n$ such that their nonincreasing rearrangement
$f^*$ satisfies $f^* (+\infty) = 0$.
\begin {theorem} [{\cite [Theorem~1] {lerner2004}}]
\label {mmsharpd}
$$
\int |f g| \leqslant c \int \left(M^\sharp_\lambda f\right) \, (M g)
$$
for any $f \in S_0$ and locally summable function $g$
with some $c$ and $\lambda$ independent of $f$ and $g$.
\end {theorem}
 %Now we are ready to state an extension of Proposition~\ref {czobounda1p}.
\begin {theorem}
\label {dpmxy}
Suppose that $X$, $Y$ and $Z$ are Banach lattices on $\mathbb R^n$ having the Fatou property,
$S_0$ is dense in $X$,
and suppose that the Hardy-Littlewood
maximal operator $M$ acts boundedly from $X$ to $Z$ and from $Y'$ to $Z'$.
Then any 
operator $T$ that satisfies estimate \eqref {flsharpczo} acts boundedly from $X$ to $Y$.
\end {theorem}
Theorem~\ref {dpmxy} 
follows at once from Theorem~\ref {mmsharpd}, since for any $f \in X \cap S_0$ and $g \in Y'$ we have
the estimate
\begin {multline}
\label {dpmxye}
\int |(T f) g| \leqslant c \int M^\sharp_\lambda (T f) \, M g \leqslant
c_1 \int M f \, M g \leqslant
\\
c_1 \| M f \|_Z \| M g \|_{Z'} \leqslant c_2 \| f \|_X \|g \|_{Y'}
\end {multline}
with some $c$, $c_1$ and $c_2$ independent of $f$ and $g$.

Since $M$ is a positive operator and $M g \geqslant g$ almost everywhere for any locally summable $g$, conditions of Theorem~\ref {dpmxy}
imply that $X \subset Z$ and $Y' \subset Z'$, %in the sense of continuous inclusions,
which in turn implies that $X \subset Z \subset Y$.
Unlike the case $X = Y = Z$ it is presently unclear whether Theorem~\ref {dpmxy} admits a converse similar to
Theorem~\ref {a1necc} below.
In other words, if a suitably nondegenerate Calderon-Zygmund operator $T$ acts boundedly from $X$ to $Y$,
does it follow that $X \subset Y$ and
there exists a lattice $Z$ satisfying the conditions of Theorem~\ref {dpmxy}?

\end {abridged}

\subsection* {Acknowledgement}

The author is grateful to S.~V.~Kisliakov who provided useful remarks to early versions of this paper
and to A.~Yu.~Kar\-lo\-vich who pointed out relevant results from \cite {karlovichlerner2005} and \cite {varlsp}.

\bibliographystyle {plain}

\bibliography {bmora}

\end {document}